\documentclass[10pt,twoside]{article}
\usepackage{graphicx}
\usepackage{amsmath,amssymb,amscd}
\usepackage{Latex-document}

\newcommand\R{\mathbb R}
\newcommand\grad{\operatorname{grad}}
\renewcommand\div{\operatorname{div}}
\newcommand\curl{\operatorname{curl}}
\newcommand\eps{\operatorname{\epsilon}}

\def\volumeno{I}

\title{\bf Differential Complexes and \vskip -2mm Numerical Stability\vskip 6mm}

\author{Douglas N. Arnold\vspace*{-0.5cm}\thanks{Institute for
Mathematics and its Applications, University of Minnesota, 400 Church
St. S.E., Minneapolis, MN 55455, USA. E-mail: arnold@ima.umn.edu}}

\date{\vspace{-8mm}}

\begin{document}
\maketitle

\thispagestyle{first} \setcounter{page}{137}

\begin{abstract}

\vskip 3mm

Differential complexes such as the de Rham
complex have recently come to play an important role in the design and
analysis of numerical methods for partial differential equations.  The
design of stable discretizations of systems of partial differential
equations often hinges on capturing subtle aspects of the structure of
the system in the discretization.  In many cases the differential
geometric structure captured by a differential complex has proven
to be a key element, and a discrete differential complex which
is appropriately related to the original complex is essential.  This
new geometric viewpoint has provided a unifying understanding of a
variety of innovative numerical methods developed over recent decades
and pointed the way to stable discretizations of problems for which
none were previously known, and it appears likely to play an important
role in attacking some currently intractable problems in numerical PDE.

\vskip 4.5mm


\end{abstract}

\vskip 12mm

\section{Introduction}\label{sc:intro}

\vskip-5mm \hspace{5mm}

During the twentieth century chain complexes, their exactness
properties, and commutative diagrams involving them pervaded many
branches of mathematics, most notably algebraic topology and
differential geometry.  Recently such homological techniques have come
to play an important role in a branch of mathematics often thought
quite distant from these, numerical analysis.   Their most significant
applications have been to the design and analysis of numerical methods
for the solution of partial differential equations.

Let us consider a general problem, such as a boundary value
problem in partial differential equations, as an operator equation:
given data $f$ in some space $Y$ find the solution $u$ in some space
$X$ to the problem $Lu=f$.  A numerical method discretizes this problem
through the construction of an operator $L_h:X_h\to Y_h$ and data
$f_h\in Y_h$ and defines an approximate solution $u_h\in X_h$ by the
equation $L_hu_h=f_h$.  Of course the numerical method is not likely to
be of value unless it is consistent which means that $L_h$ and $f_h$
should be close to $L$ and $f$ in an appropriate sense.

Before we speak of solving the original problem, numerically or
otherwise, we should first confront the question of whether it is
well-posed.  That is, given $f\in Y$, does a unique $u\in X$ exist,
and, if so, do small changes of $f$ induce small changes in $Y$?  The
analogous questions for the numerical method, whether given $f_h\in Y_h$
a unique $u_h\in X_h$ is determined by the discrete equation
$L_hu_h=f_h$, and whether small changes in $f_h$ induce small changes
in $u_h$, is the question of \emph{stability} of the numerical method.
A common paradigm, which can be formalized in many contexts of
numerical analysis, is that a method which is consistent and stable is
convergent.

Well-posedness is a central issue in the theory of partial differential
equations.  Of course, we do not expect just any PDE problem to be
well-posed. Well-posedness hinges on structure of the problem which may
be elusive or delicate.  Superficially small changes, for example to
the sign of a coefficient or the type of boundary conditions, can
certainly destroy well-posedness.  The same is true for the stability
of numerical methods: it often depends on subtle or elusive properties
of the numerical scheme.  Usually stability reflects some portion of
the structure of the original problem that is captured by the numerical
scheme.  However in many contexts it is not enough that the numerical
scheme be close to the original problem in a quantitative sense for it
to inherit stability.  That is, it may well happen that a consistent
method for a well-posed problem is unstable.  In this paper we shall
see several examples where the exactness properties of discrete
differential complexes and their relation to differential complexes
associated with the PDE are crucial tools in establishing the
stability of numerical methods.  In some cases the homological
arguments have served to elucidate or validate methods that had been
developed over the preceding decades.  In others they have pointed the
way to stable discretizations of problems for which none were
previously known.  They will very likely play a similar role in the
eventual solution of some formidable open problems in numerical PDE,
especially for problems with significant geometric content, such as in
numerical general relativity. As in other branches of mathematics, in
numerical analysis differential complexes serve both to encode key
structure concisely and to unify considerations from seemingly very
different contexts.

In this paper we shall discuss only finite element methods since, of
the major classes of numerical methods for PDE, they are the most
amenable to rigorous analysis, and have seen the greatest use of
differential complexes.  But complexes have recently arisen in the study of
finite differences, finite volumes, and spectral methods as well.

\section{Finite element spaces}\label{sc:fem}

\vskip-5mm \hspace{5mm}

A finite element space on a domain $\Omega$ is a function space
defined piecewise by a certain assembly procedure which we now
recall; cf.~\cite{ciarlet}. For simplicity, here we shall restrict
to spaces of piecewise polynomials with respect to a triangulation
of an $n$-dimensional domain by $n$-simplices with $n=2$ or $3$
(so implicitly we are assuming that $\Omega\subset\R^2$ is
polygonal or $\Omega\subset\R^3$ is polyhedral). On each simplex
$T$ we require that there be given a function space of \emph{shape
function} $W_T$ and a set of \emph{degrees of freedom}, i.e., a
set of linear functionals on $W_T$ which form a basis for the dual
space.  Moreover, each degree of freedom is supposed to be
associated with some subsimplex of some dimension, i.e., in three
dimensions with a vertex, an edge, a face, or the tetrahedron
itself.  For a subsimplex which is shared by two simplices in the
triangulation, we assume that the corresponding functionals are in
one-to-one-correspondence.  Then the finite element space $W_h$ is
defined as those functions on $\Omega$ whose restriction to each
simplex $T$ of the triangulation belongs to $W_T$ and for which
the corresponding degrees of freedom agree whenever a subsimplex
is shared by two simplices.

The simplest example is obtained by choosing $W_T$ to be the constant
functions and taking as the only degree of freedom on $T$ the $0$th
order moment
$\phi\mapsto \int_T\phi(x)\,dx$ (which we associate with $T$ itself).
The resulting finite element space is simply the space of piecewise
constant functions with respect to the given triangulation.  Similarly
we could choose $W_T=\mathbb P_1(T)$ (by $\mathbb P_p(T)$ we denote the
space of polynomial functions on $T$ of degree at most $p$), and take
as degrees of freedom  the  moments of degrees $0$ and also those of
degree $1$,
$\phi\mapsto\int_T\phi(x)x_i\,dx$.  Again all the degrees of freedom
are associated to $T$ itself.  This time the finite element space
consists of all piecewise linear functions.  Of course, the
construction extends to higher degrees.

\begin{figure}[ht!]
\centerline{\includegraphics[height=1.25in]{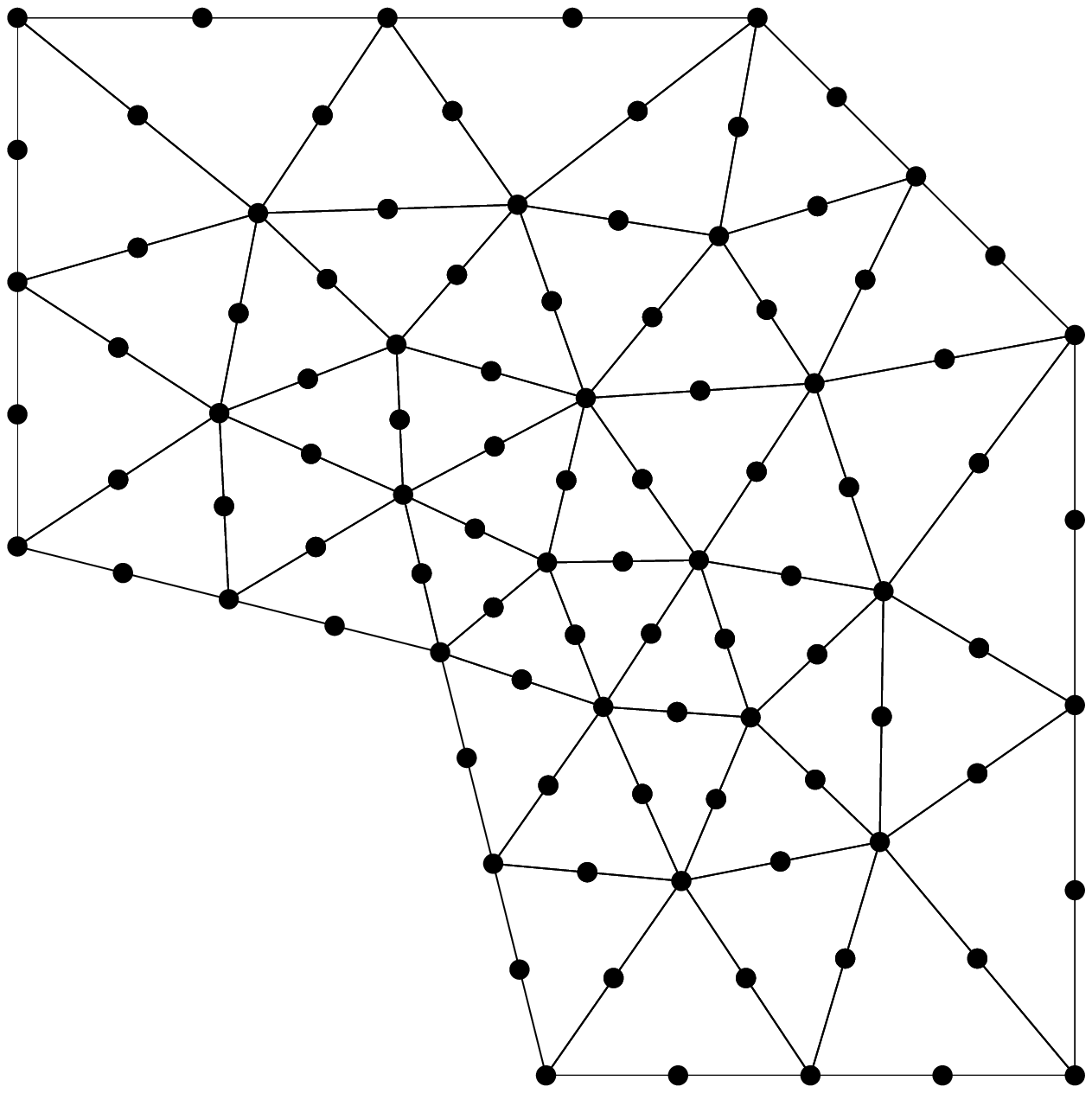}\quad
\includegraphics[height=1.25in]{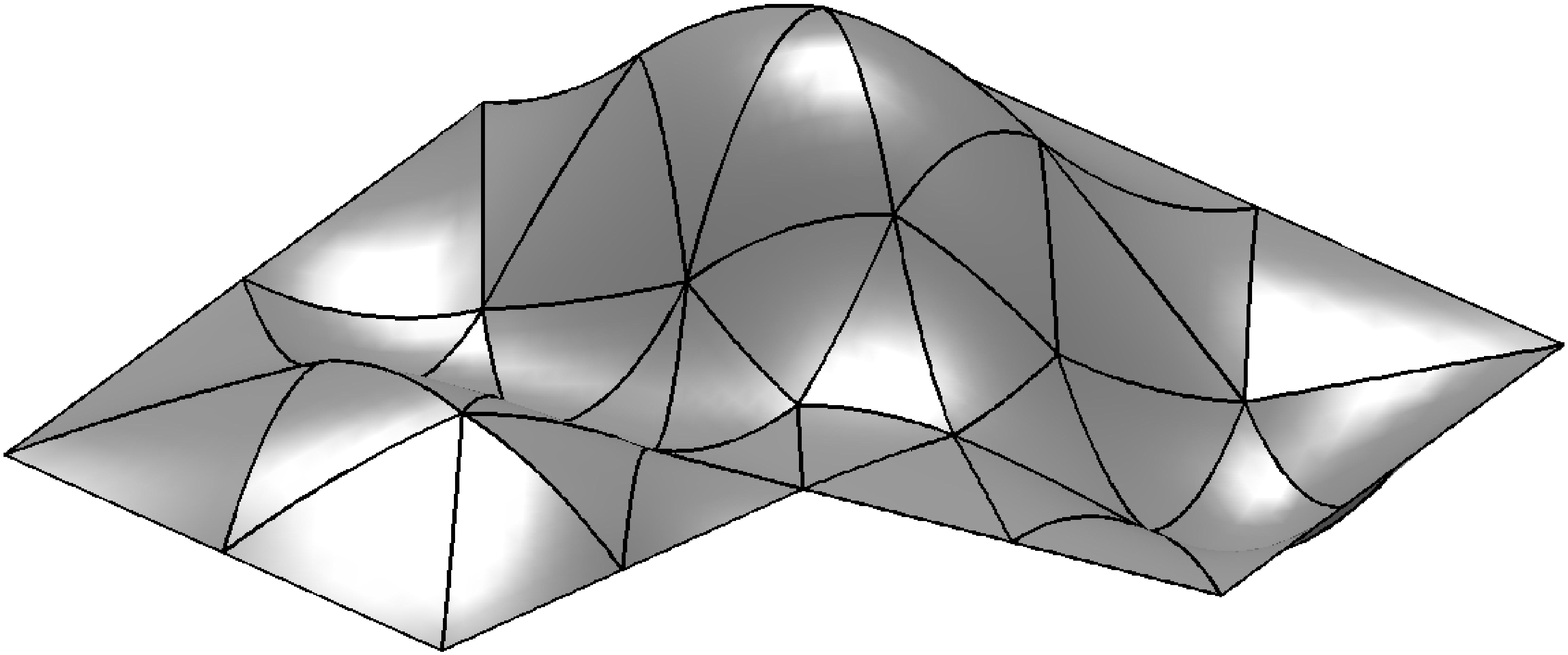}}
\begin{center}
\begin{minipage}{10cm}
\caption{\label{fg:assemb} A mesh marked with the locations of the
degrees of freedom for Lagrange finite elements of degree $2$ and
a typical such finite element function.}
\end{minipage}
\end{center}
\end{figure}

A more common piecewise linear finite element space occurs if we again
choose $W_T=\mathbb P_1(T)$, but take as degrees of freedom the maps
$\phi\mapsto \phi(v)$, one associated to each vertex $v$.  In this case
the assembled finite element space consists of all \emph{continuous}
piecewise linear functions.  More generally we can choose $W_T=\mathbb
P_p(T)$ for $p\ge 1$, and associate to each vertex the evaluation
degrees of freedom just mentioned, to each edge the moments on the edge
of degree at most $p-2$, to each face the moments on the face of degree
at most $p-3$, and to each tetrahedron the moments of degree at most
$p-4$.  The resulting finite element space, called the \emph{Lagrange
finite element} of degree $p$, consists of all continuous piecewise
polynomials of degree at most $p$.  Figure~\ref{fg:assemb} shows a mesh
of a two dimensional domain and a typical function in the space
of Lagrange finite elements of degree $2$ with respect to this mesh.

\begin{figure}[ht!]
\centerline{\includegraphics[width=1in]{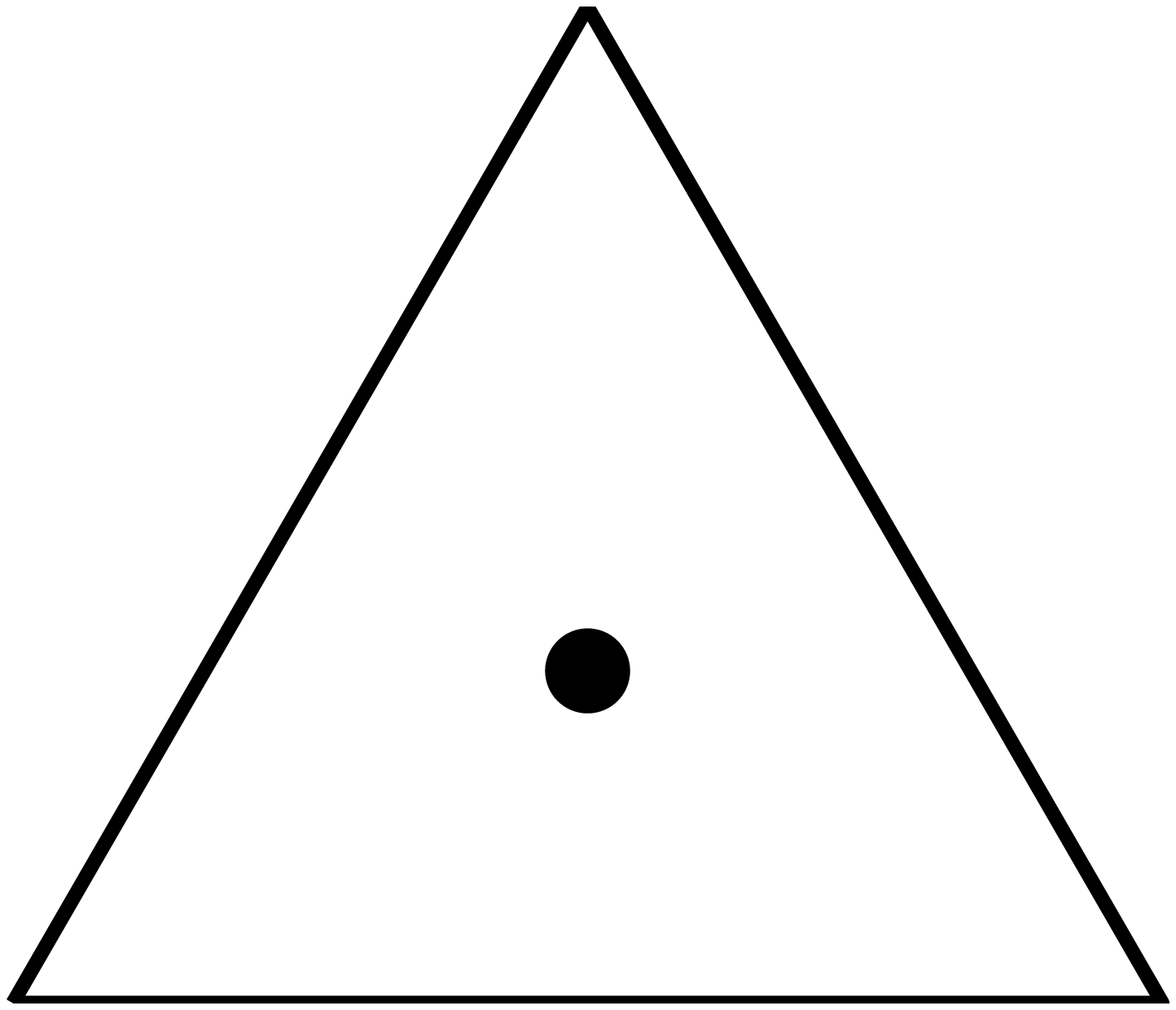}\quad
\includegraphics[width=1in]{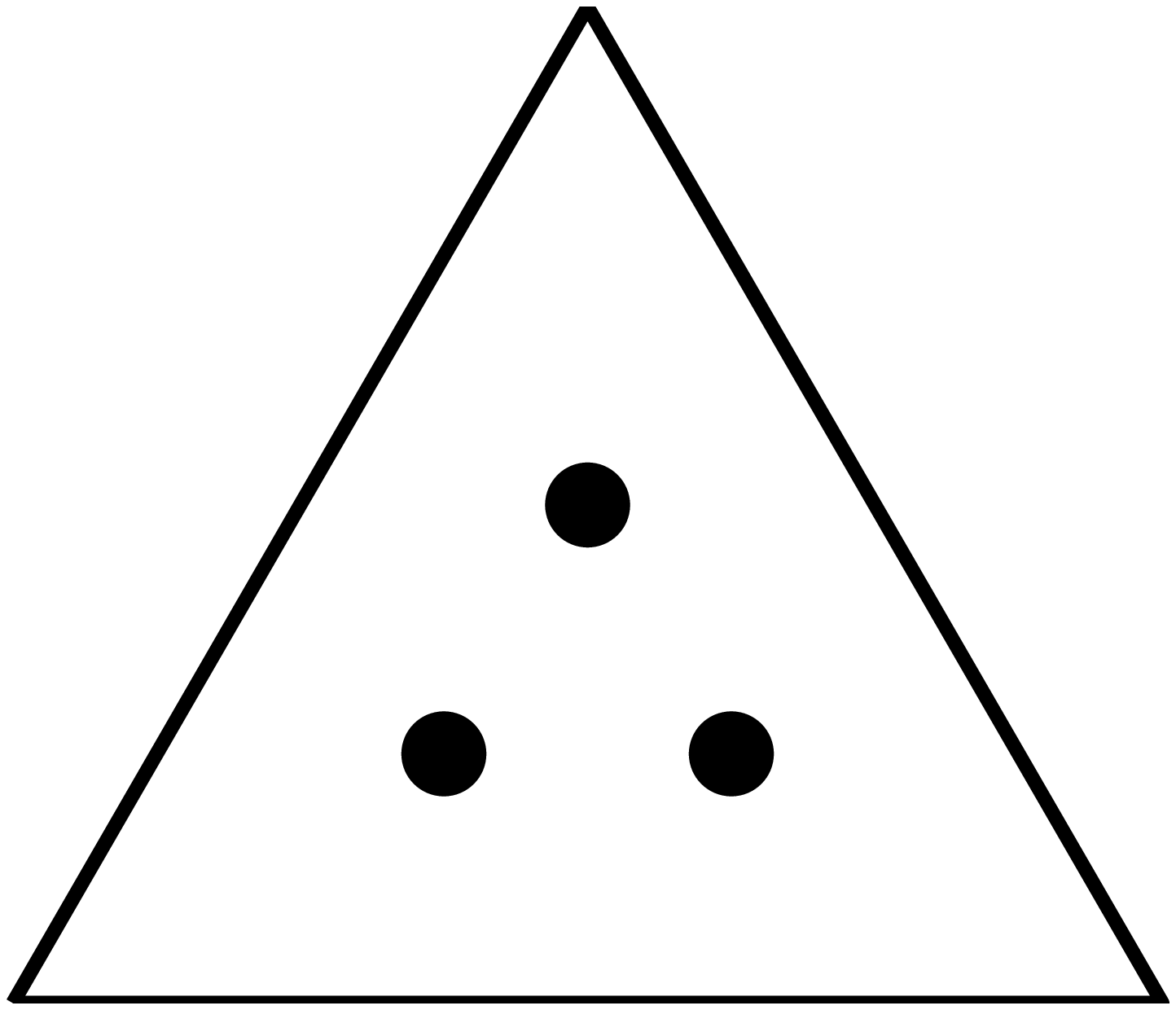}\quad
\includegraphics[width=1in]{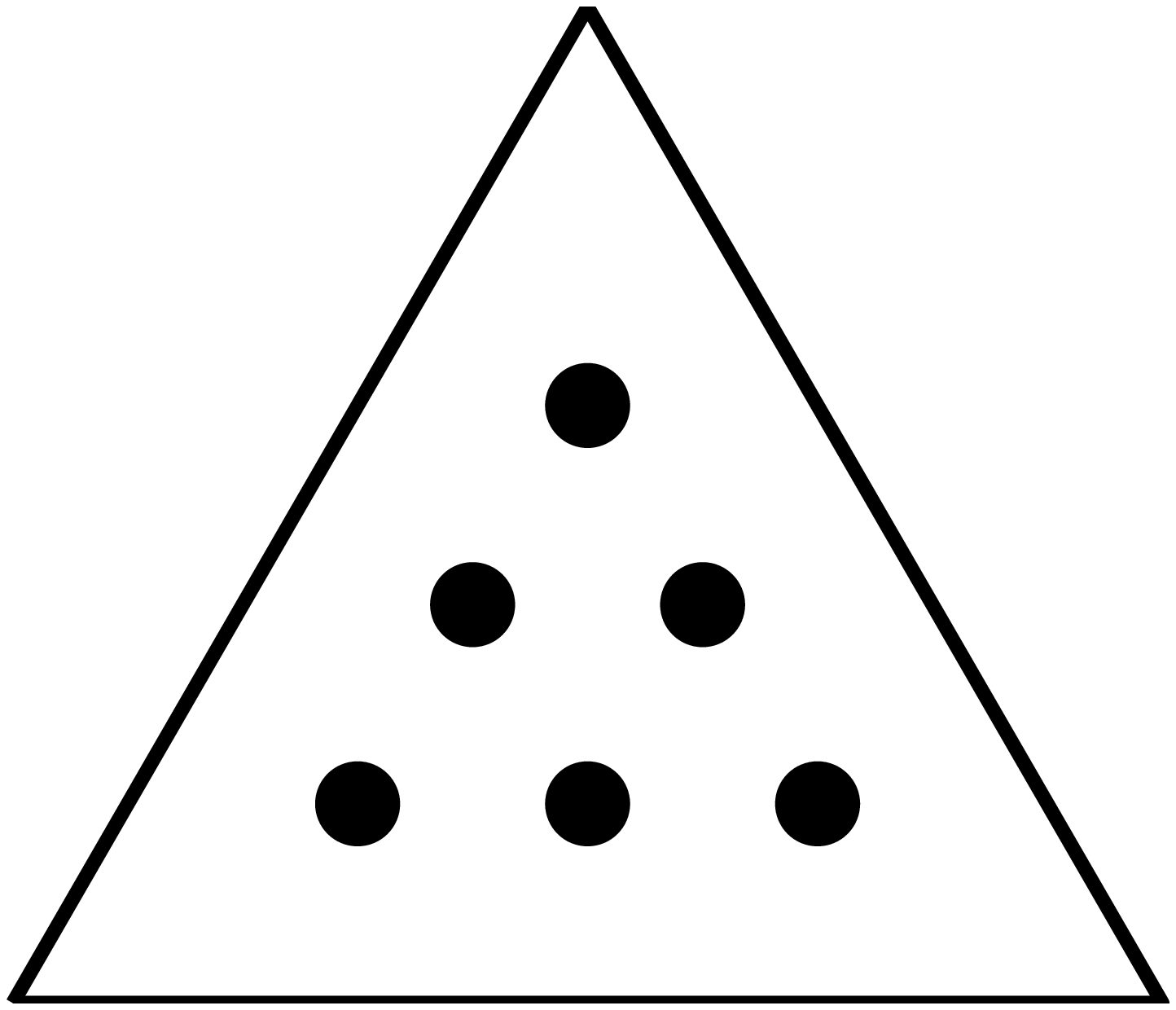}}
\vskip.1in \centerline{\includegraphics[width=1in]{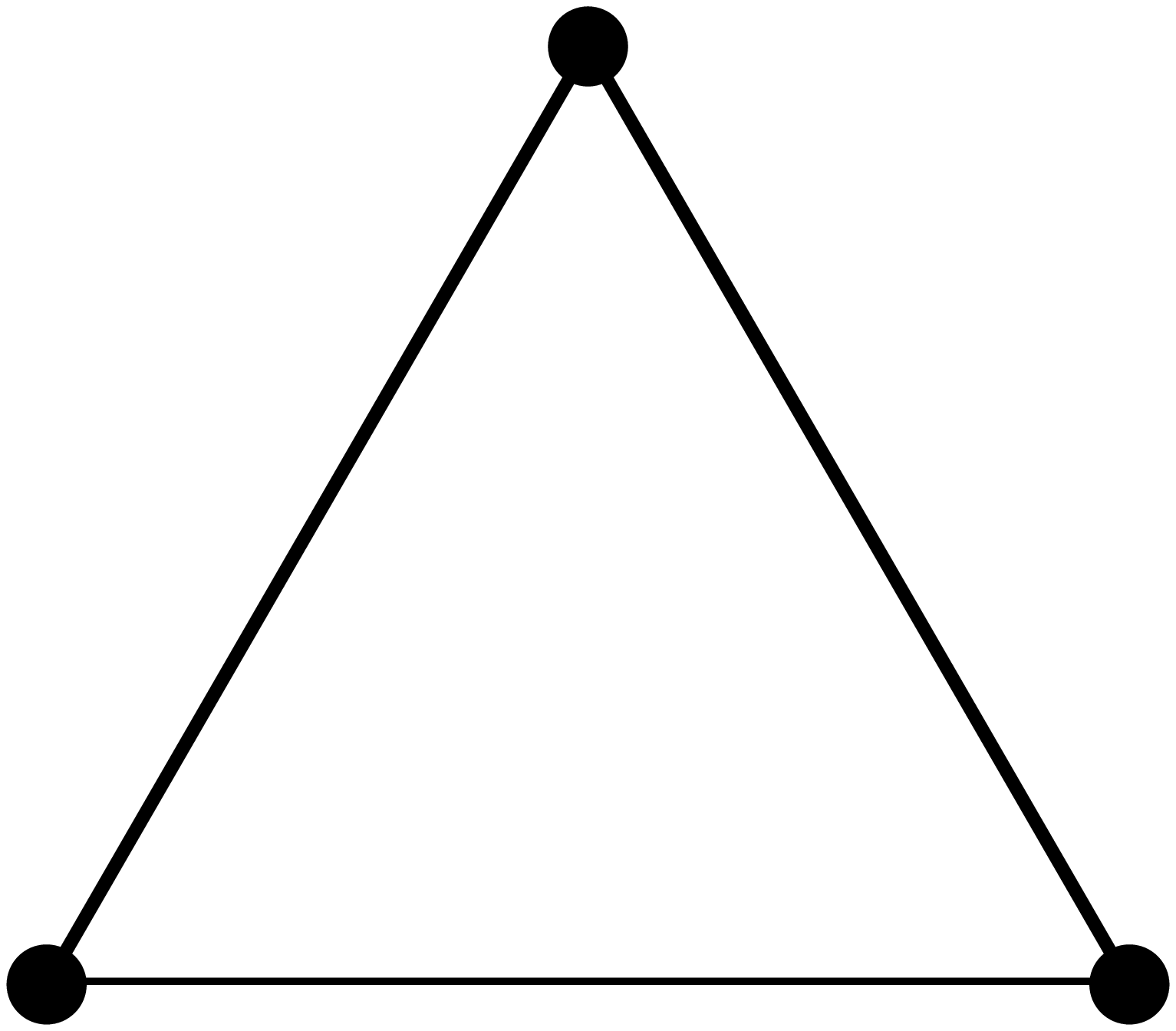}\quad
\includegraphics[width=1in]{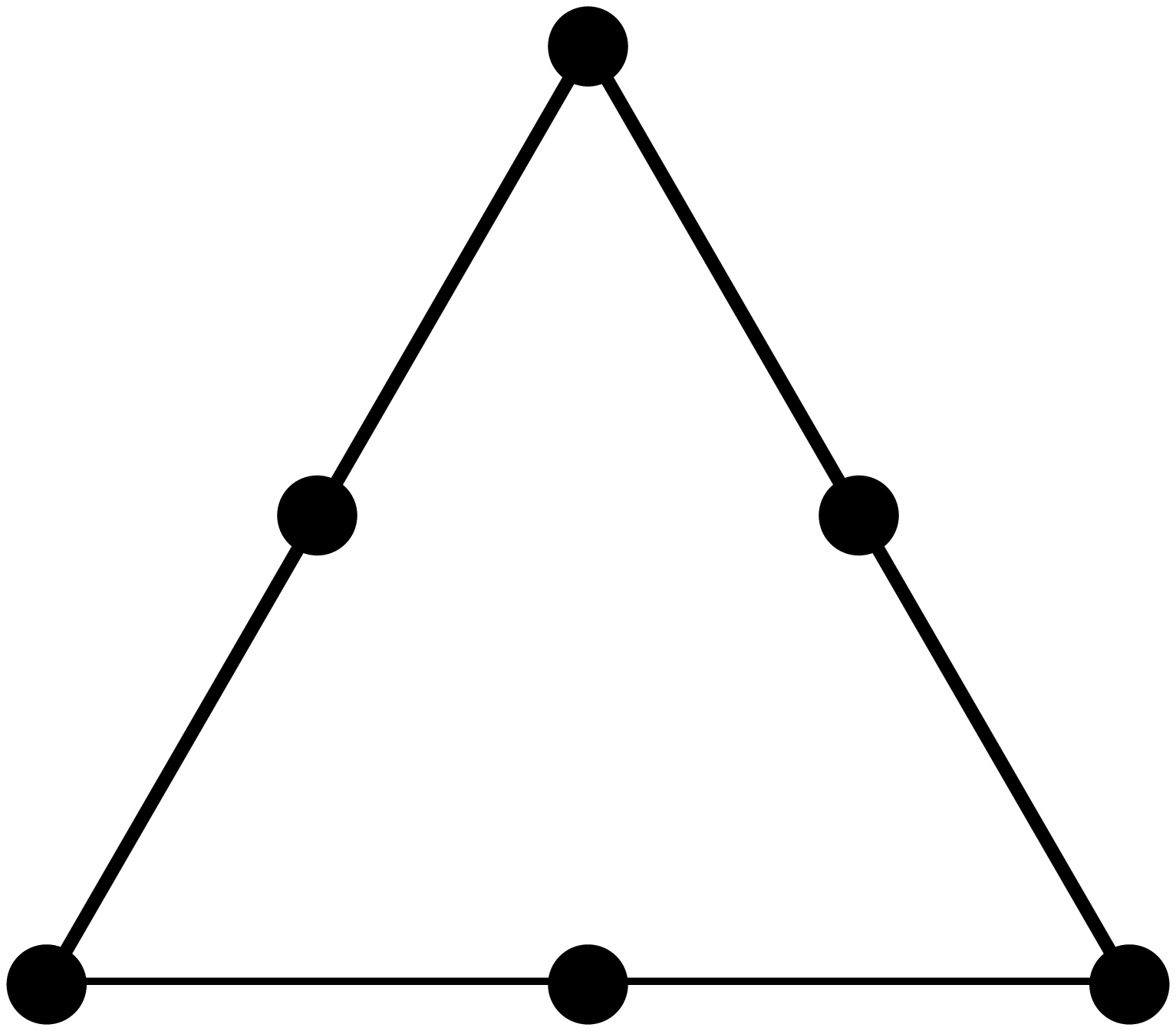}\quad
\includegraphics[width=1in]{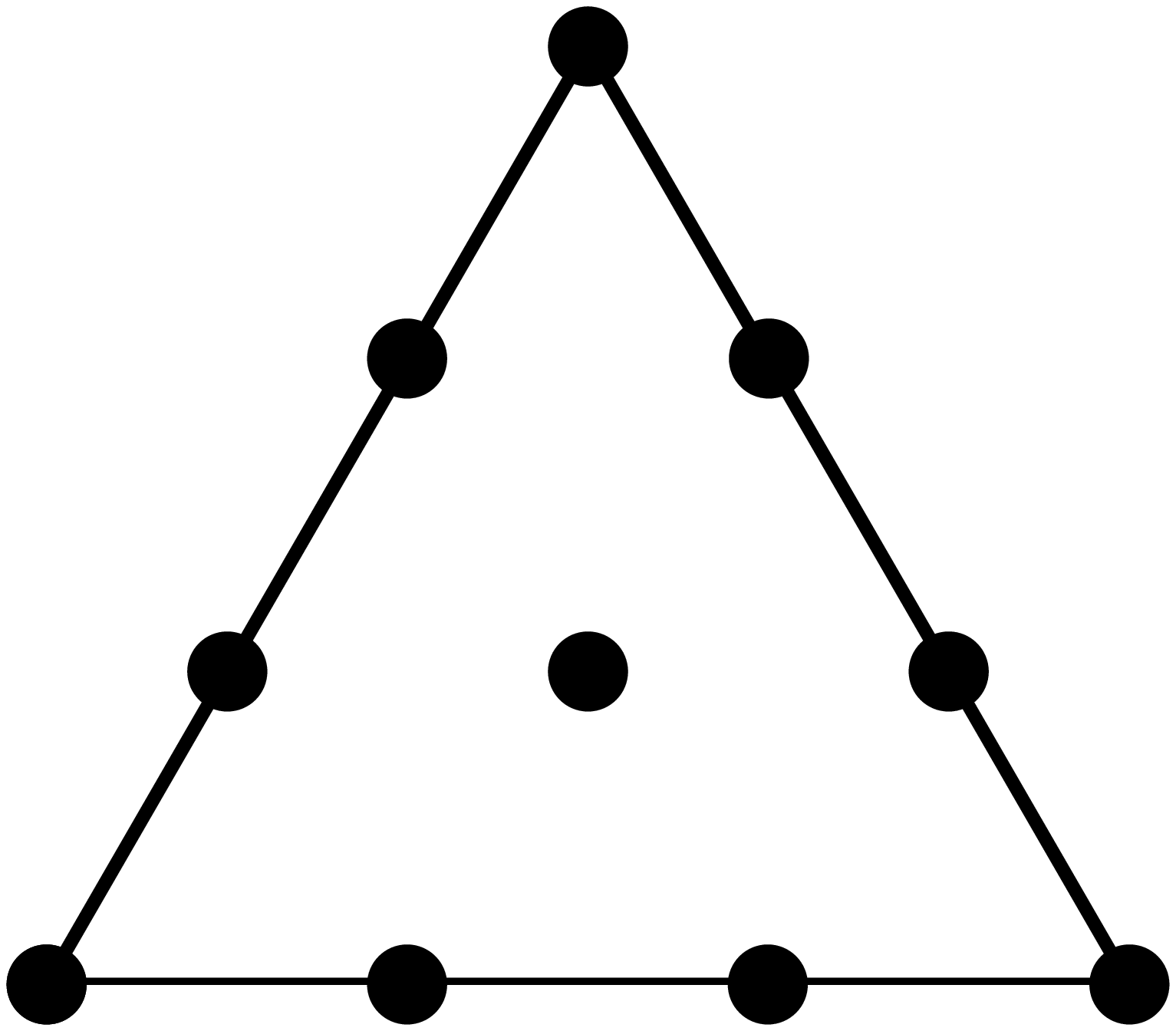}}
\vskip.1in \centerline{\includegraphics[width=1in]{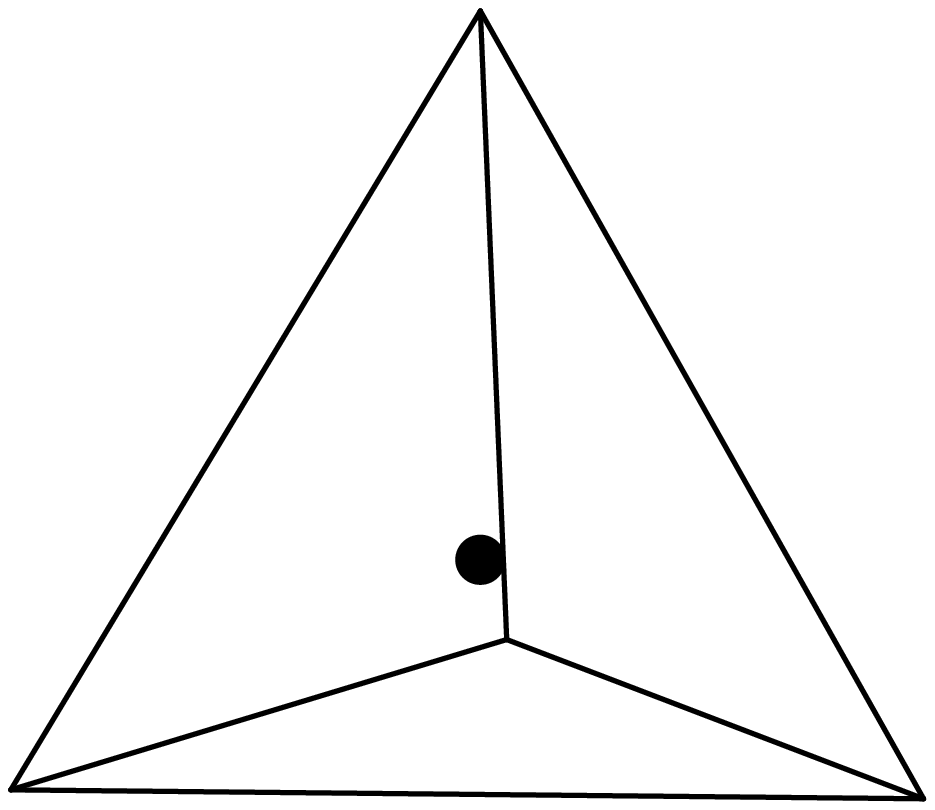}\quad
\includegraphics[width=1in]{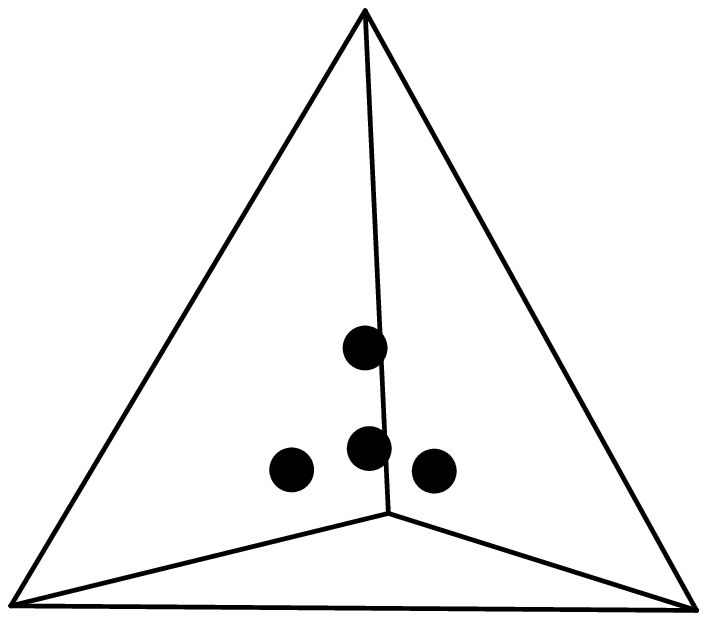}\quad
\includegraphics[width=1in]{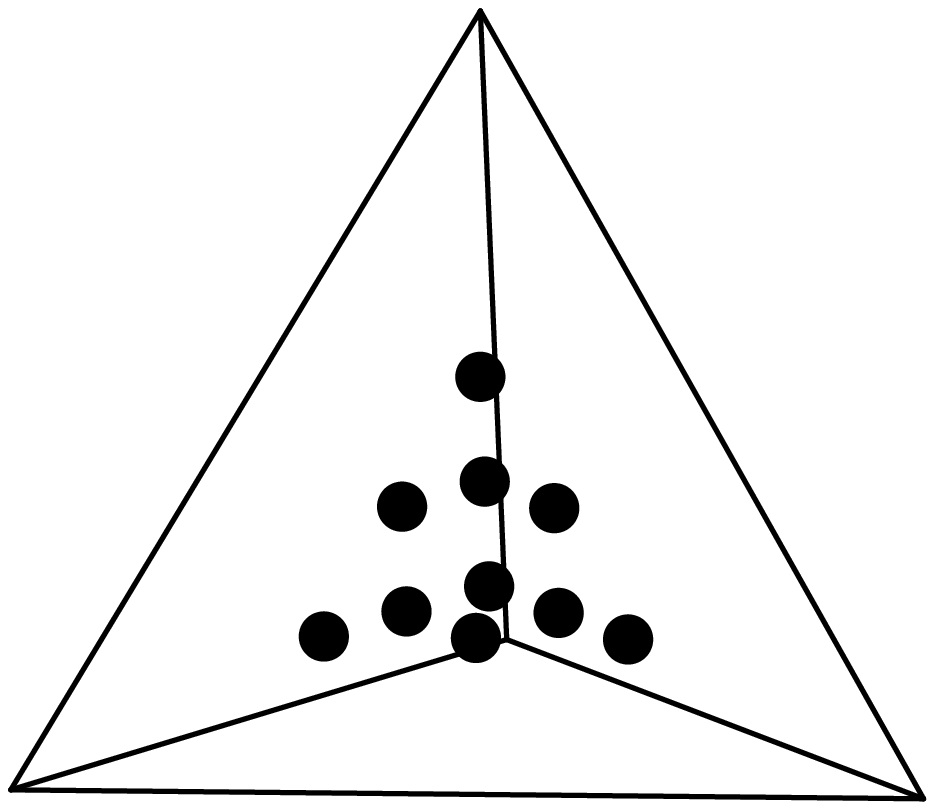}}
\vskip.1in \centerline{\includegraphics[width=1in]{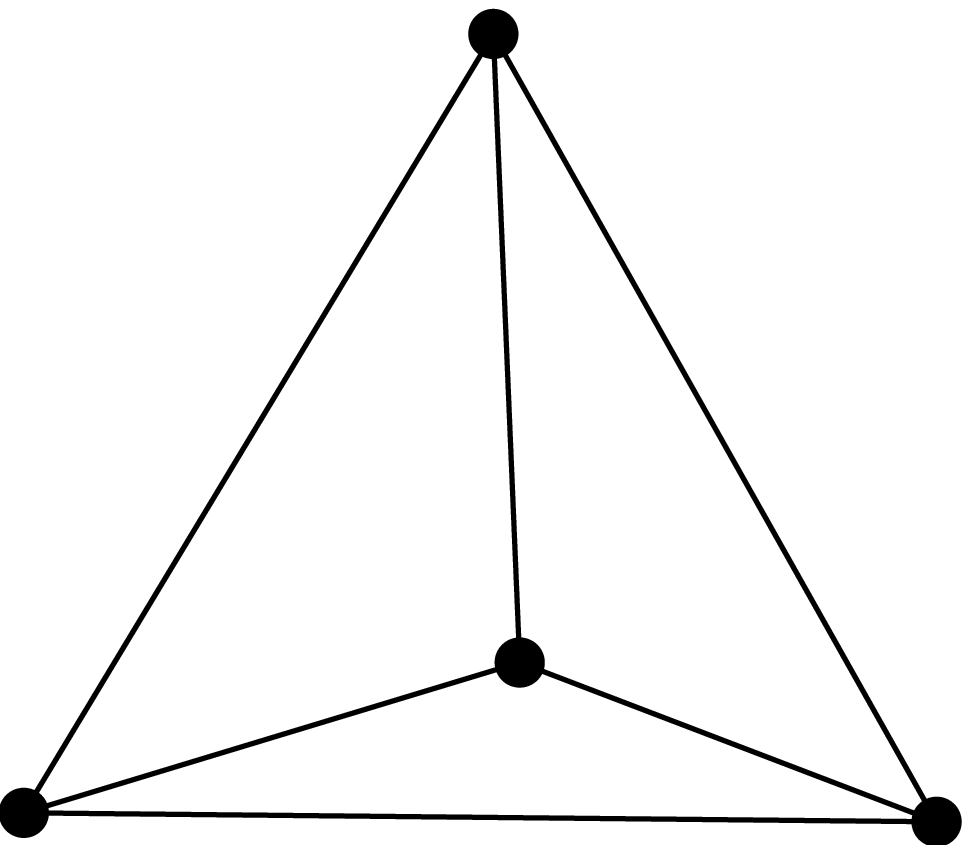}\quad
\includegraphics[width=1in]{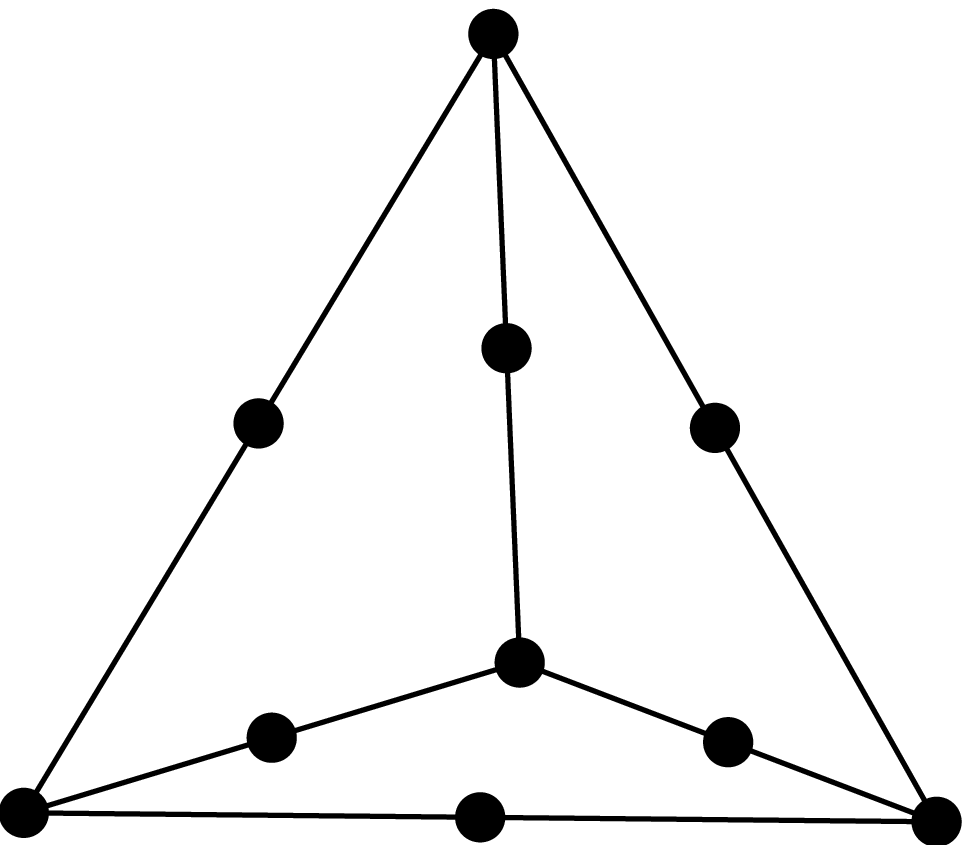}\quad
\includegraphics[width=1in]{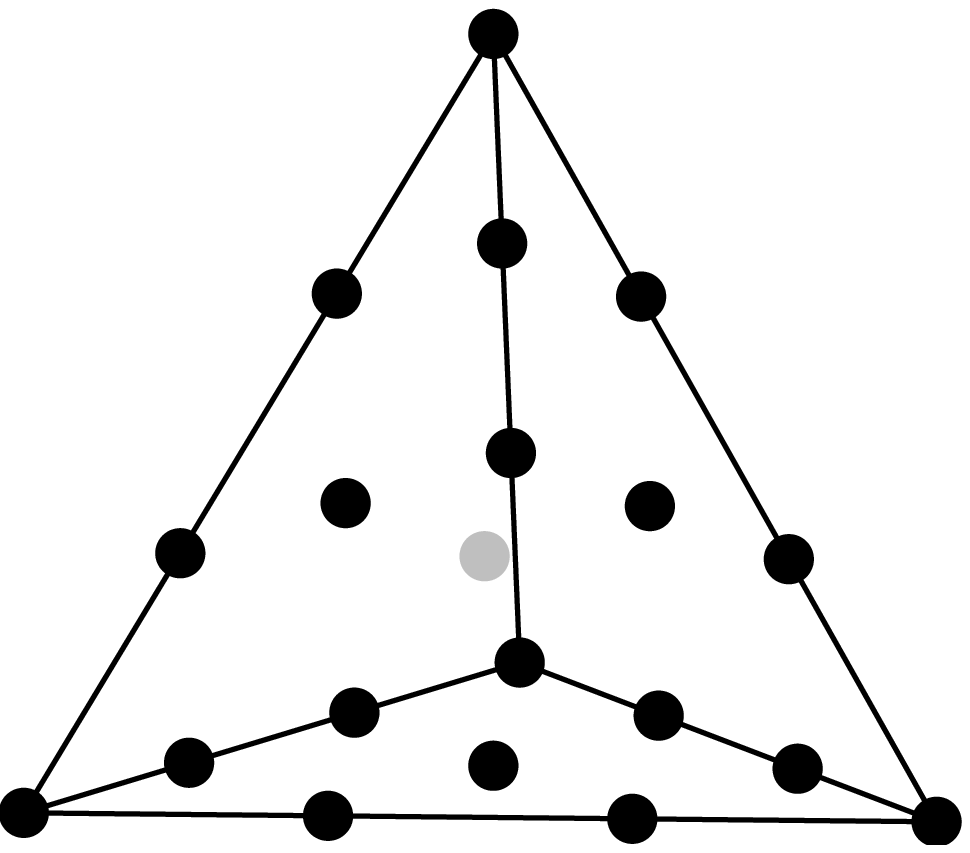}}
\begin{center}
\begin{minipage}{10cm}
\caption{\label{fg:elts1} Element diagrams.  First row:
discontinuous elements of degrees $0$, $1$, and $2$ in two
dimensions. Second row: Lagrange elements of degrees $1$, $2$, and
$3$ in two dimensions. Third and fourth rows: the corresponding
elements in three dimensions.}
\end{minipage}
\end{center}
\end{figure}

Mnemonic diagrams as in Figure~\ref{fg:elts1} are often associated to
finite element spaces, depicting a single element $T$ and a marker for
each degree of freedom.

\begin{figure}[ht!]
\centerline{\includegraphics[width=1in]{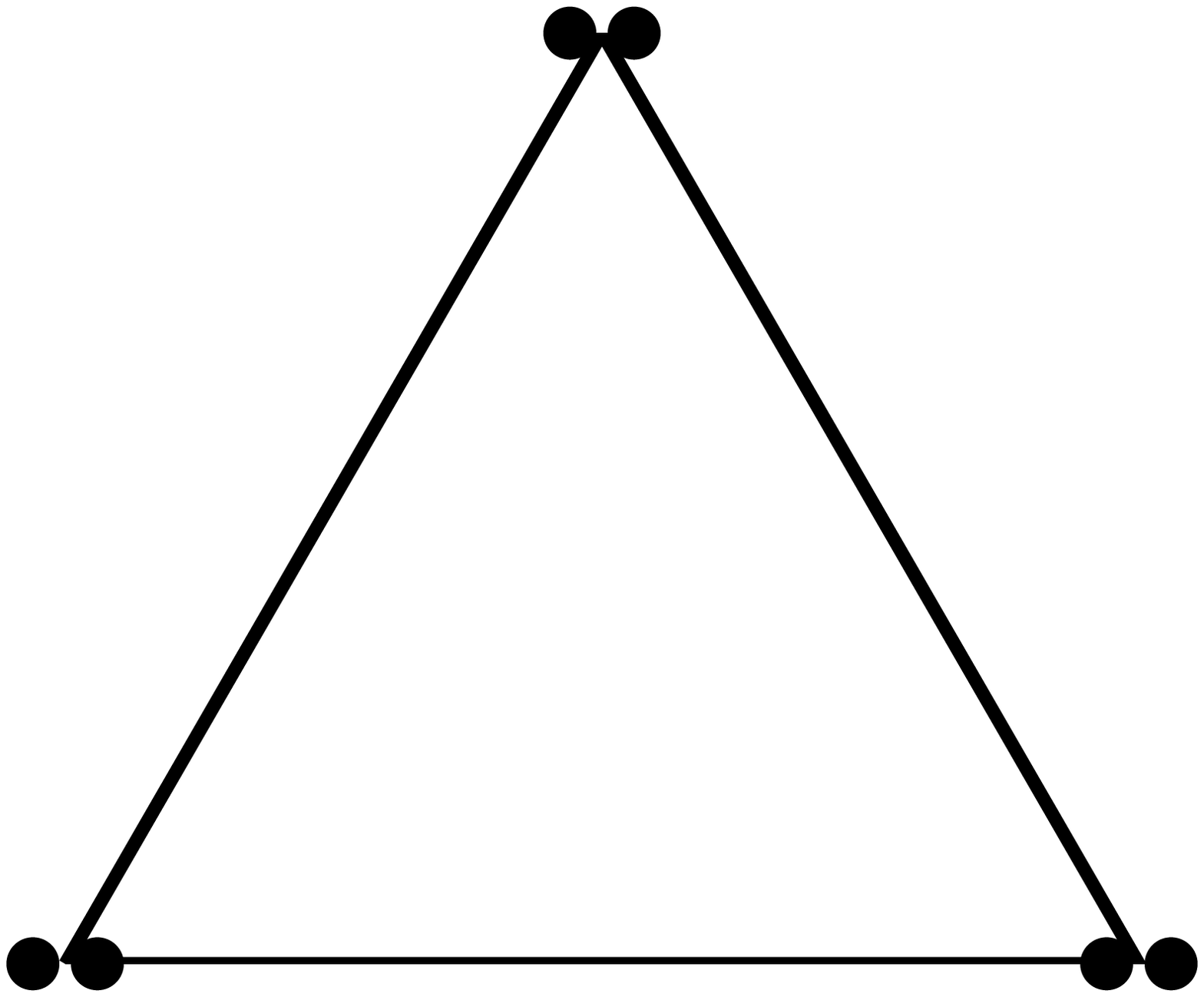}\quad
\raise-.18in\hbox{\includegraphics[width=1in]{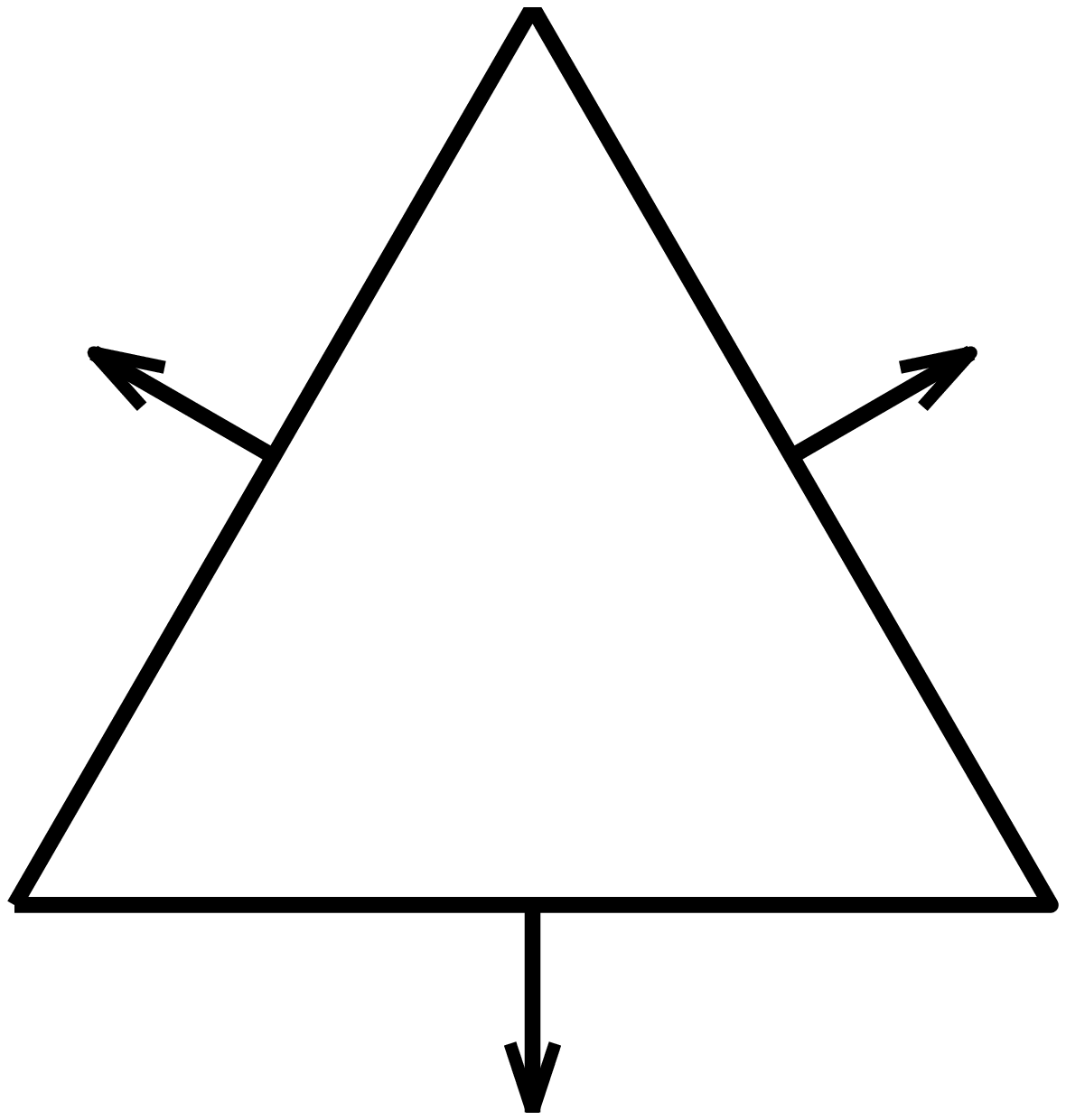}}\quad
\includegraphics[width=1in]{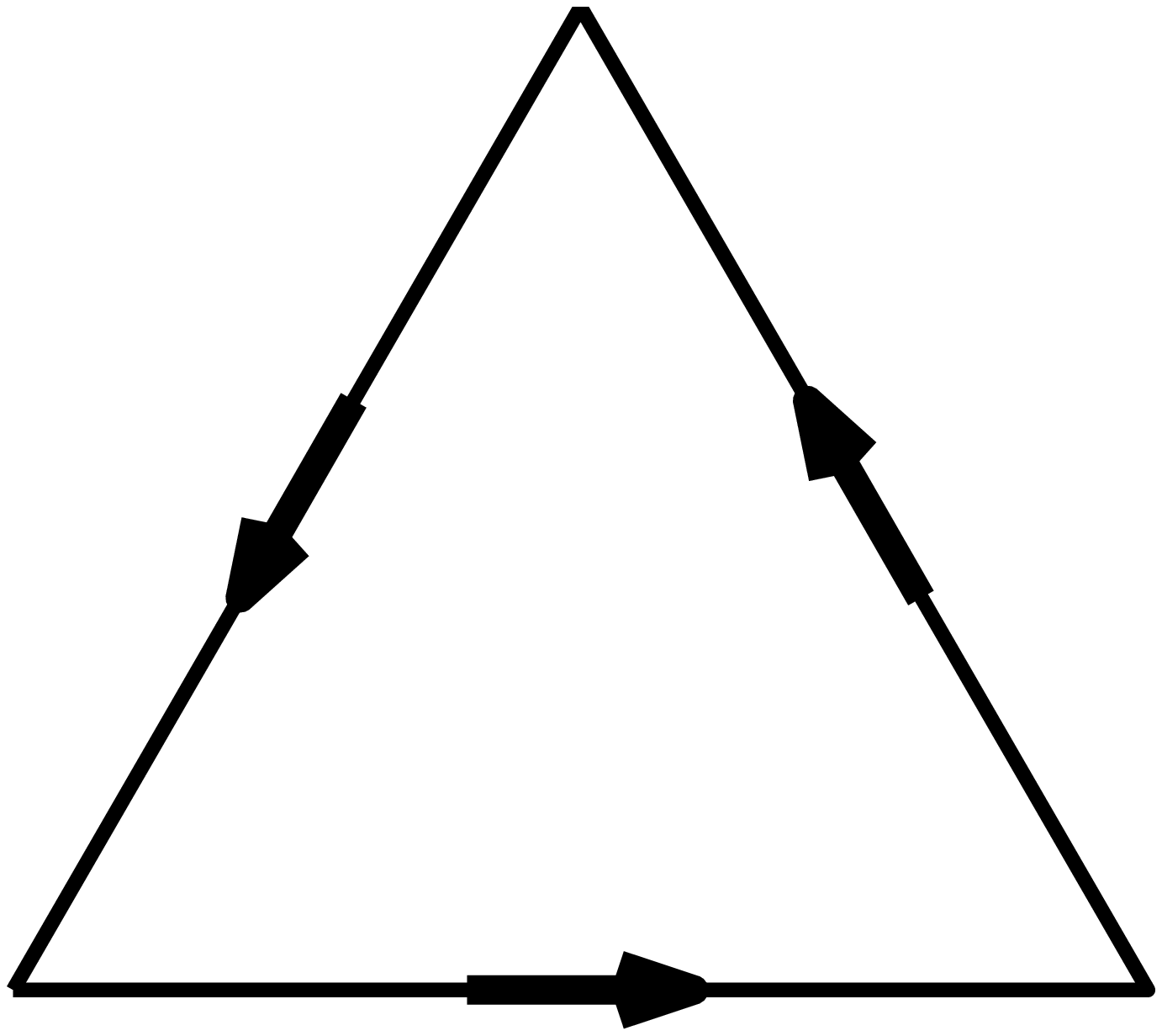}}
\vskip.1in \centerline{\includegraphics[width=1in]{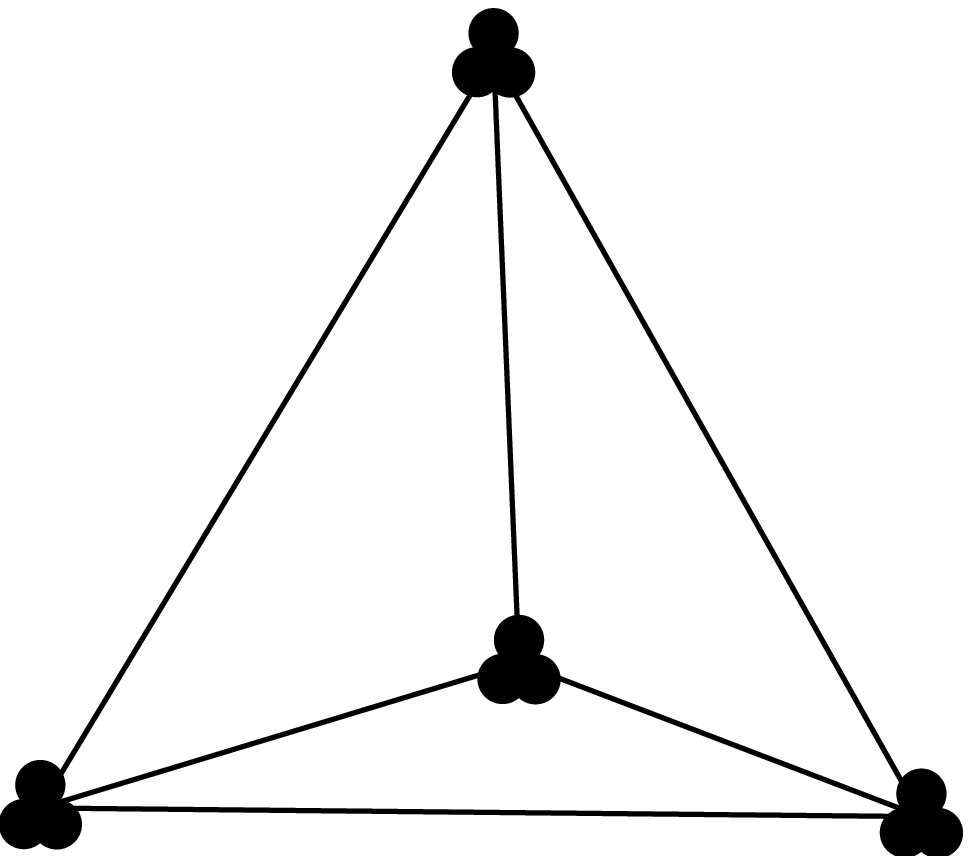}\quad
\raise-.2in\hbox{\includegraphics[width=1in]{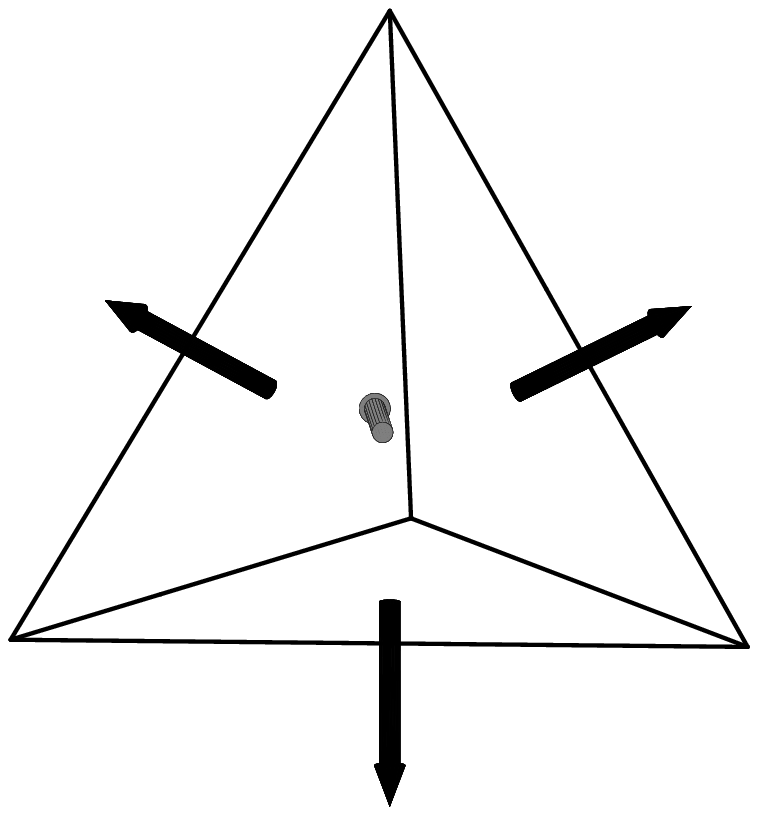}}\quad
\includegraphics[width=1in]{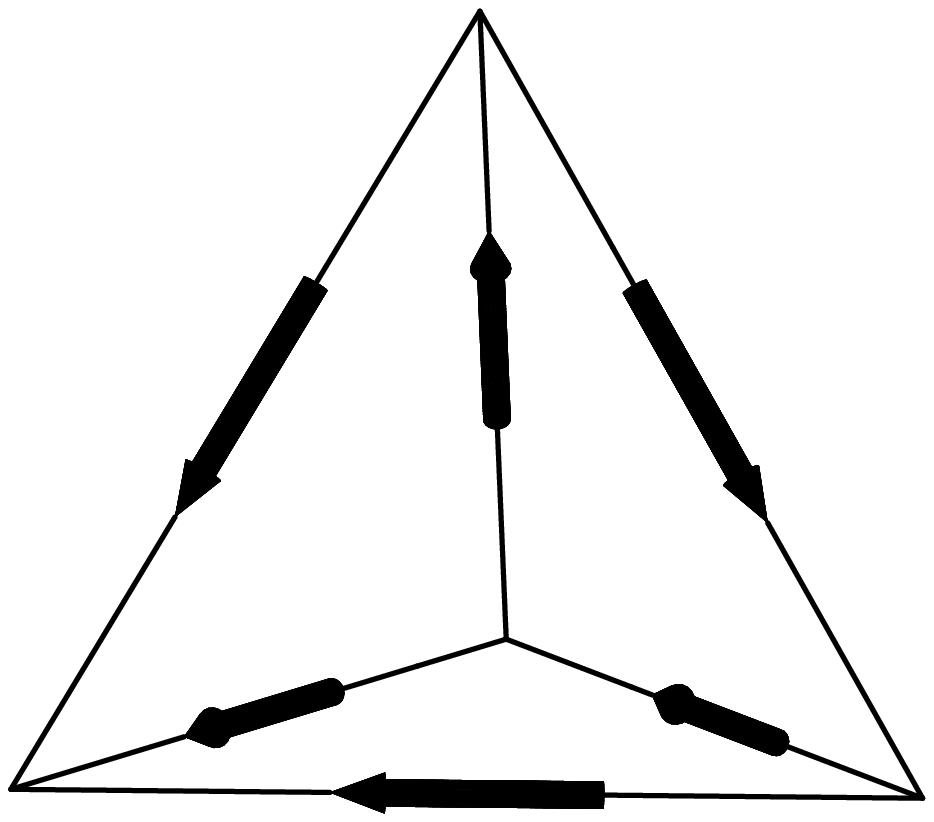}}
\begin{center}
\begin{minipage}{10cm}
\caption{\label{fg:eltsvec} Element diagrams for some finite
element approximations to vector fields in two and three
dimensions. Multiple dots are used as markers to indicate the
evaluation of all components of a vector field.  Arrows are used
for normal moments on codimension $1$ subsimplices and for
tangential components on edges. Left: continuous piecewise linear
fields.  Middle: face elements of lowest order.  Right: edge
elements of lowest order.}
\end{minipage}
\end{center}
\end{figure}

Next we describe some finite element spaces that can be used to
approximate vector-valued functions.  For brevity we limit the
descriptions to the $3$-dimensional case, but supply diagrams in both
$2$ and $3$ dimensions. Of course we may simply take the Cartesian
product of three copies of one of the previous spaces.  For example,
the element diagrams shown on the left of Figure~\ref{fg:eltsvec} refer
to continuous piecewise linear vector fields in two and three
dimensions.  More interesting spaces are the \emph{face elements} and
\emph{edge elements} essentially conceived by Raviart and Thomas
\cite{raviart-thomas} in two dimensions and by Nedelec \cite{nedelec}
in three dimensions.
In the lowest order case, the face elements take as shape functions
polynomial vector fields of the form $p(x)=a+bx$ where $a\in\R^3$,
$b\in\R$ and $x=(x_1,x_2,x_3)$, a 4-dimensional subspace of the
12-dimensional space  $\mathbb P_1(T,\R^3)$ of polynomial vector fields
of degree at most $1$.  The degrees of freedom are taken to be the $0$th
order moments of the
normal components on the faces of codimension $1$,  $p\mapsto \int_f
p(x)\cdot n_f\,dx$ where $f$ is a face and $n_f$ the unit normal to the
face. The element diagram is shown in the middle column of
Figure~\ref{fg:eltsvec}. In the lowest order case the edge elements
shape functions are polynomial vector fields of the form $p(x)=a+b\times
x$ where $a,b\in\R^3$, which form a 6-dimensional subspace of  $\mathbb P_1(T,\R^3)$.  The degrees of
freedom are the $0$th order moments over the edges of the component tangent to
the edge, $p\mapsto \int_e p(x)\cdot t_e\,dx$, as indicated on the
right of Figure~\ref{fg:eltsvec}.

\begin{figure}[ht!]
\centerline{\raise-.16in\hbox{\includegraphics[width=1in]{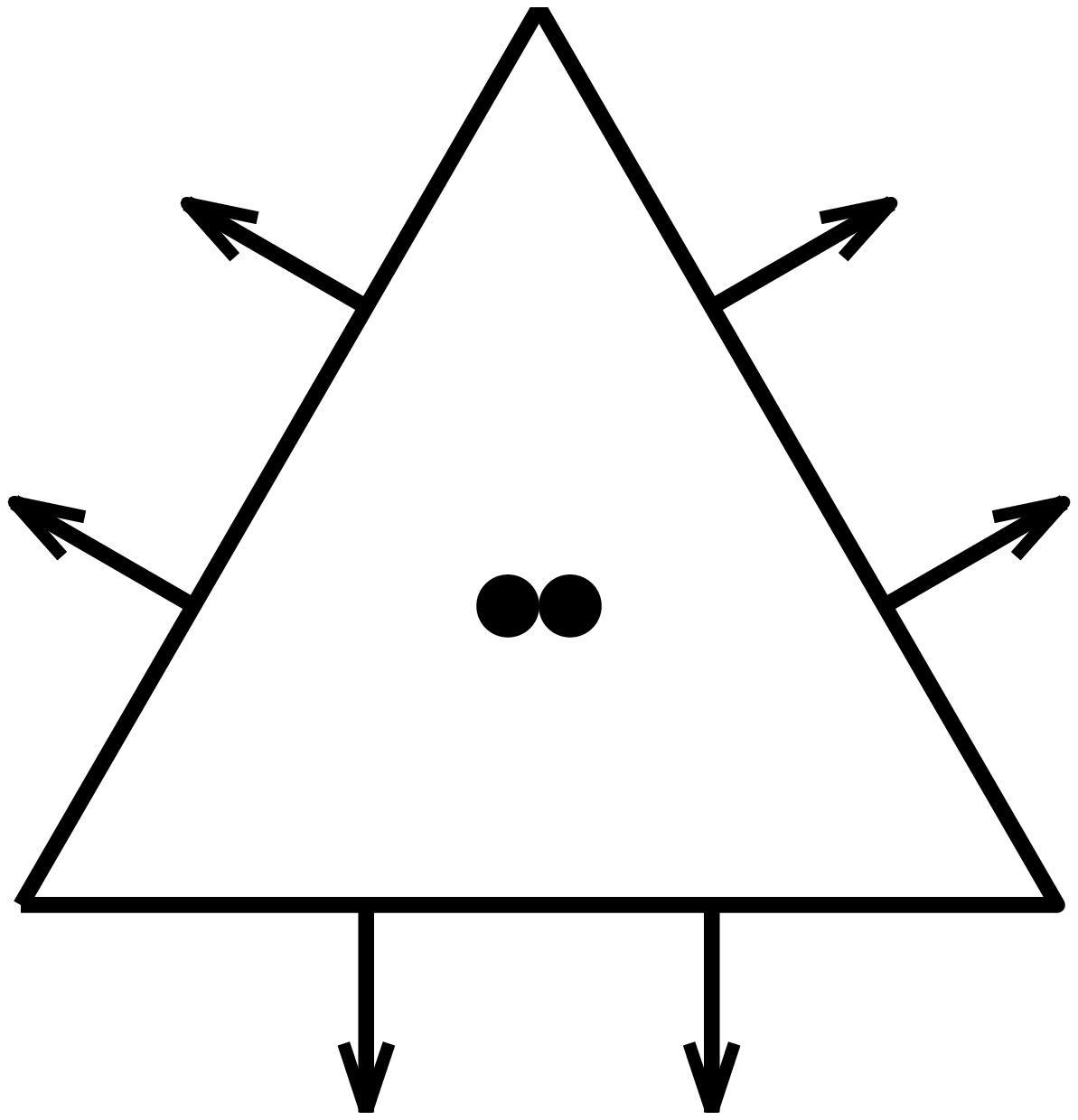}}\quad
\raise-.20in\hbox{\includegraphics[width=1in]{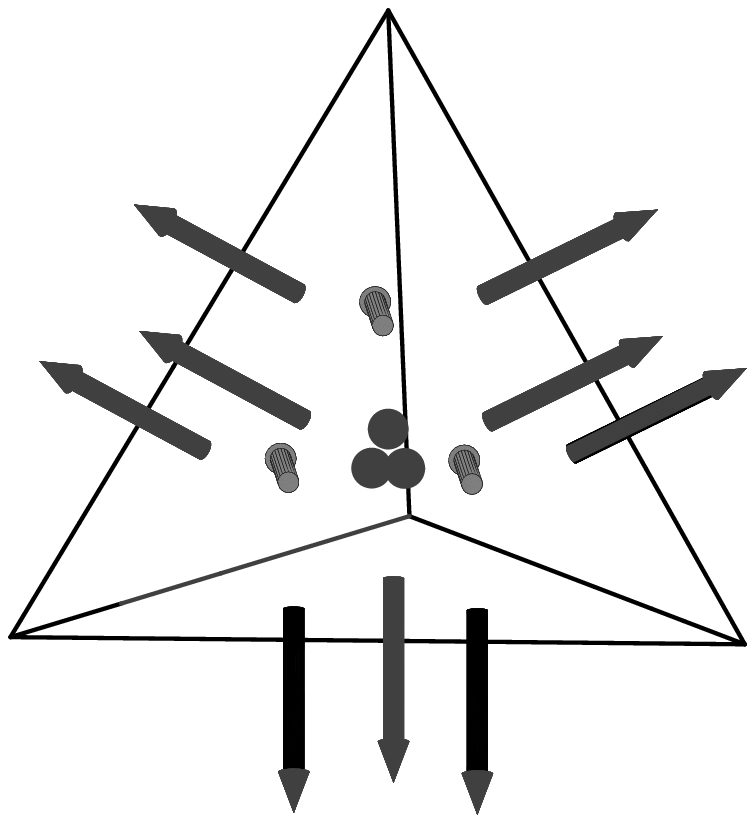}}\quad
\raise-.01in\hbox{\includegraphics[width=1in]{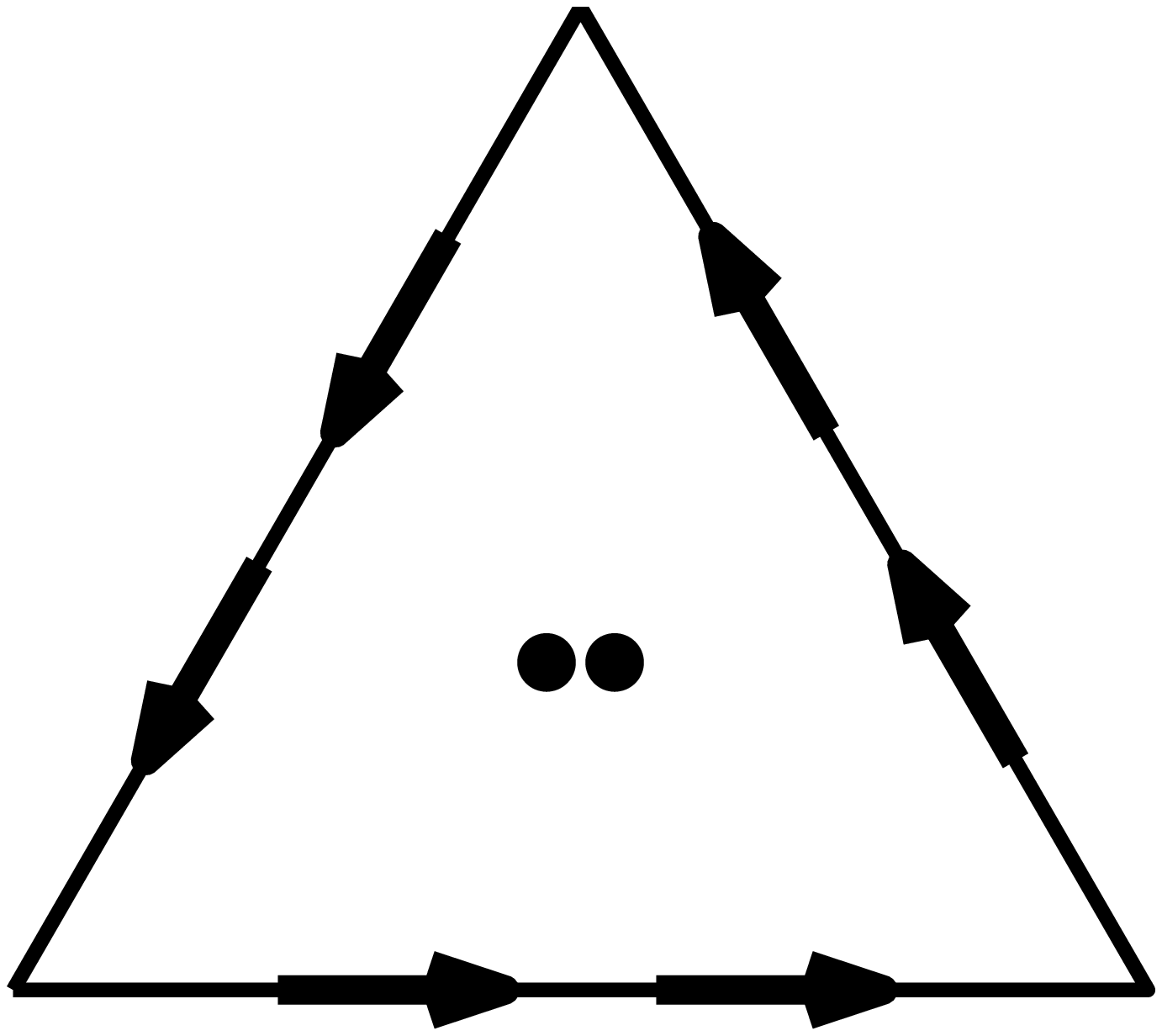}}\quad
\raise-.0in\hbox{\includegraphics[width=1in]{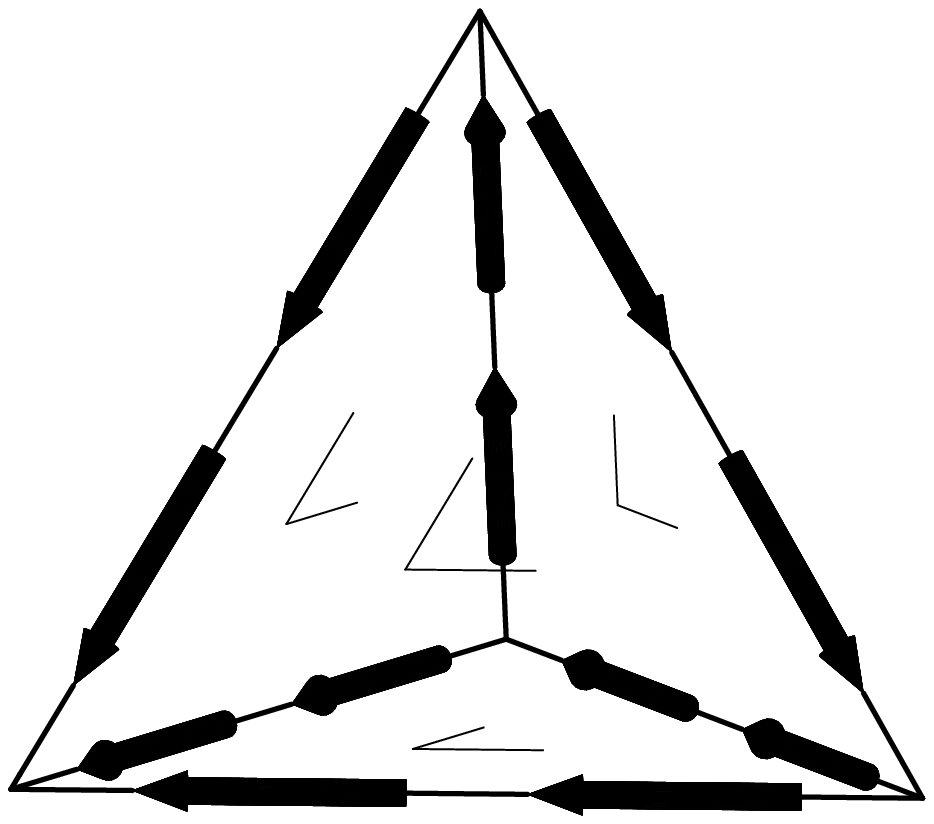}}}
\begin{center}
\begin{minipage}{10cm}
\caption{\label{fg:edgeface} The face (left) and edge (right)
elements of the second lowest order in 2- and 3-dimensions.}
\end{minipage}
\end{center}
\end{figure}

Each of these spaces can be generalized to arbitrarily high order. For
the next higher order face space, the shape functions take the form
$p(x)=a(x)+b(x)x$ where $a\in\mathbb P_1(T,\R^3)$ and
$b\in\mathbb P_1(T)$ a linear scalar-valued polynomial.  This gives a subspace of
$\mathbb P_2(T,\R^3)$ of dimension $15$,
and the degrees of freedom are the moments of degree at most $1$ of the
normal components on the faces and the moments of degree $0$ of all
components on the tetrahedron. This element is indicated on the left of
Figure~\ref{fg:edgeface}. For the second lowest order edge space, the shape
functions take the form $p(x)=a(x)+b(x)\times x$ with
$a,b\in \mathbb P_1(T,\R^3)$, giving a $20$-dimensional space.
The degrees of freedom are the tangential moments of degree at most $1$
on the edges (two per edge) and the tangential moments of degree $0$ on
the faces (two per face).  This element is indicated on the right of
Figure~\ref{fg:edgeface}.

The choice of the shape functions and the degrees of freedom determine the
smoothness of the functions belonging to the assembled finite element
space.  For example, the Lagrange finite element spaces
of any degree belong to the Sobolev space $H^1(\Omega)$ of $L^2(\Omega)$
functions whose distributial first partial derivatives also belong
to $L^2(\Omega)$ (and even to
$L^infty(\Omega)$).  In fact, the distributional first partial
derivative of a continuous piecewise smooth function coincides with its
derivative taken piecewise and so belongs to $L^2$.
Thus the degrees of freedom we imposed in constructing the Lagrange
finite elements are sufficient to insure that the assembled finite
element space $W_h\subset H^1(\Omega)$.  In fact more is true:
for the Lagrange finite element space with shape function spaces
$W_T=\mathbb P_p(T)$, we have
\begin{equation*}
W_h = \{\, u\in H^1(\Omega)\,|\,u|_T\in W_T \text{ for all simplices
$T$ of the triangulation}\,\}.
\end{equation*}
This says that, in a sense, the degrees of freedom impose exactly
the continuity required to belong to $H^1$, no less and no more.

In contrast, the discontinuous piecewise polynomial spaces are subsets
of $L^2(\Omega)$ but not of $H^1(\Omega)$, since their distributional
first derivatives involve distributions supported on the interelement
boundaries, and so do not belong to $L^2(\Omega)$.

For the vector-valued finite elements there are more possibilities. The
face and edge spaces contain discontinuous functions, and so are not
contained in $H^1(\Omega,\R^3)$.  However, for vector fields belonging
to one of the face spaces the normal component of the vector field does
not jump across interelement boundaries, and this implies, via
integration by parts, that the distributional divergence of the
function coincides with the divergence taken piecewise.  Thus the face
spaces belong to $H(\div,\Omega)$, the space of $L^2$ vector fields on
$\Omega$ whose divergence belongs to $L^2$.  Indeed, for these spaces
the degrees of freedom impose exactly the continuity of $H(\div)$, no
less or more.  For the edge spaces it can be shown that the tangential
components of a vector field do not jump across element boundaries,
and this implies that the edge functions belong to $H(\curl,\Omega)$,
the space of $L^2$ vector fields whose curl belongs to $L^2$.  Again
the degrees of freedom impose exactly the continuity needed for
inclusion in $H(\curl)$.

\section{Discrete differential complexes}

\vskip-5mm \hspace{5mm}

The de~Rham complex
\begin{equation*}
\begin{CD}
\R  \hookrightarrow\, \bigwedge^0(\Omega) @>d>>
 \bigwedge^1(\Omega) @>d>> \cdots @>d>> \bigwedge^n(\Omega) \to0
\end{CD}
\end{equation*}
is defined for an arbitrary smooth $n$-manifold $\Omega$.  Here
$\bigwedge^k(\Omega)$ denotes the space of differential $k$-forms
on $\Omega$, i.e., for $\omega\in \bigwedge^k(\Omega)$ and
$x\in\Omega$, $\omega(x)$ is an alternating $k$-linear map on the
tangent space $T_x\Omega$. The operators
$d:\bigwedge^k(\Omega)\to\bigwedge^{k+1}(\Omega)$ denote exterior
differentiation. This is is a complex in that the composition of
two exterior differentiations always vanishes. Moreover, and if
the manifold is topologically trivial, then it is exact.

If $\Omega$ is a domain in $\R^3$, then we may identify its  tangent
space at any point with $\R^3$.  Using the Euclidean inner product, the
space of linear maps on $\R^3$ may be identified by $\R^3$ as usual, so
$\bigwedge^1(\Omega)$ may be identified with the space
$C^\infty(\Omega,\R^3)$ of smooth
vector fields on $\Omega$.  Moreover, the space of alternating bilinear
maps on $\R^3$ may be identified with $\R^3$ by associating to
a vector $u$ the alternating bilinear map $(v,w)\mapsto \det(u|v|w)$.
Thus we have an identification of $\bigwedge^2(\Omega)$ with $\R^3$ as
well.  Finally the only alternating trilinear maps on $\R^3$ are
given by multiples of the determinant map $(u,v,w)\mapsto c\det(u|v|w)$,
and so we may identify $\bigwedge^3(\Omega)$ with $C^\infty(\Omega)$.  In
terms of such \emph{proxy fields}, the de~Rham complex becomes
\begin{equation}\label{eq:derham}
\begin{CD}
\R\hookrightarrow\, C^\infty(\Omega) @>\grad>>
 C^\infty(\Omega,\R^3) @>\curl>> C^\infty(\Omega,\R^3) @>\div>>
 C^\infty(\Omega,\R) \to0.
\end{CD}
\end{equation}
Alternatively we may consider $L^2$-based forms and the sequence becomes
\begin{equation*}
\begin{CD}
\R\hookrightarrow\, H^1(\Omega) @>\grad>>
 H(\curl,\Omega) @>\curl>> H(\div,\Omega) @>\div>>
 L^2(\Omega,\R) \to0.
\end{CD}
\end{equation*}

The finite element spaces constructed above allow us to
form discrete analogues of the de~Rham complex. Given
some triangulation of $\Omega\subset\R^3$,
let $W_h$ denote the space of continuous piecewise linear
finite elements, $Q_h$ the lowest order edge element space,
$S_h$ the lowest order face element space, and $V_h$ the
space of piecewise constants.  Then $\grad W_h\subset Q_h$
(since $Q_h$ contains all piecewise constant vector fields
belonging to $H(\curl)$ and the gradient of a continuous
piecewise linear is certainly such a function), $\curl Q_h\subset S_h$
(since $S_h$ contains all piecewise constant vector fields belonging
to $H(\curl)$), and $\div S_h\subset V_h$.  Thus we have the discrete
differential complex
\begin{equation}\label{eq:ddc}
\begin{CD}
\R\hookrightarrow\, W_h @>\grad>> Q_h @>\curl>>  S_h @>\div>> V_h \to0.
\end{CD}
\end{equation}
This differential complex captures the topology of the domain to
the same extent as the de~Rham complex.  In particular, if the
domain is topologically trivial, then the sequence is exact.

It is convenient to abbreviate the above statement using the element
diagrams introduced earlier.  Thus we will say that the following
complex is exact:
\begin{equation*}
\begin{CD}
\R  \hookrightarrow\,@.
  \raise-.15in\hbox{\includegraphics[width=.5in]{000212.eps}} @>\grad>>
  \raise-.15in\hbox{\includegraphics[width=.5in]{000220.eps}} @>\curl>>
  \raise-.25in\hbox{\includegraphics[width=.5in]{000219.eps}} @>\div>>
  \raise-.15in\hbox{\includegraphics[width=.5in]{000209.eps}} @.\to0
\end{CD}
\end{equation*}
By this we mean that if we assemble finite element spaces $W_h$, $Q_h$,
$S_h$, and $V_h$ using the indicated finite elements and a triangulation
of a topologically trivial domain, then the corresponding discrete
differential complex \eqref{eq:ddc} is exact.

There is another important relationship between the de~Rham complex
\eqref{eq:derham} and the discrete complex \eqref{eq:ddc}.  The defining
degrees
of freedom determine projections $\Pi^W_h:C^\infty(\Omega)\to W_h$,
$\Pi^Q_h:C^\infty(\Omega,\R^3)\to Q_h$, and so on.  In fact $\Pi^W_h$ is
just the usual interpolant, $\Pi^V_h$ is the $L^2$-projection into
the piecewise constants, and the projections $\Pi^Q_h$ and
$\Pi^S_h$ onto the edge and face elements are determined by
the maintenance of the appropriate moments.  It can be checked, based on
Stokes theorem, that the following diagram commutes.
\begin{equation}\label{eq:cd}
\begin{CD}
\R  \hookrightarrow\,@. C^\infty(\Omega,\R) @>\grad>>
 C^\infty(\Omega,\R^3) @>\curl>> C^\infty(\Omega,\R^3) @>\div>>
  C^\infty(\Omega,\R)@. \to 0\\
 @.@VV\Pi^W_hV @VV\Pi^Q_hV  @VV\Pi^S_hV @VV{\Pi^V_h}V\\
\R  \hookrightarrow\,@.W_h @>\grad>>  Q_h @>\curl>>S_h @>\div>>
V_h@.\to 0
\end{CD}
\end{equation}

The finite element spaces appearing in this diagram, with one
degree of freedom for each vertex for $W_h$, for each edge for $Q_h$,
for each face for $S_h$, and for each simplex for $V_h$,
are highly geometrical.  In fact, recalling the identifications
between fields and differential forms, we may view these spaces as
spaces of piecewise smooth differential forms.  They were in fact
first constructed in this context, without any thought of finite
elements or numerical methods, by Whitney \cite{whitney}.  The spaces
were reinvented, one-by-one, as finite element spaces in response
to the needs of various numerical problems, and the properties which
are summarized in the commutative diagram above were slowly rediscovered
as needed to analyze the resulting numerical methods.  The connection
between low order edge and face finite elements and Whitney forms
was first realized by Bossavit \cite{bossavit}.

Analogous statements hold for higher order Lagrange, edge, face,
and discontinuous finite elements.  For example, the following
diagram commutes and has exact rows:
\begin{equation*}
\begin{CD}
\R  \hookrightarrow\,@. C^\infty(\Omega,\R) @>\grad>>
 C^\infty(\Omega,\R^3) @>\curl>> C^\infty(\Omega,\R^3) @>\div>>
  C^\infty(\Omega,\R)@. \to 0\\
 @.@VVV @VVV @VVV @VVV\\
\R  \hookrightarrow\,@.
  \raise-.15in\hbox{\includegraphics[width=.5in]{000213.eps}} @>\grad>>
  \raise-.15in\hbox{\includegraphics[width=.5in]{000224.eps}} @>\curl>>
  \raise-.25in\hbox{\includegraphics[width=.5in]{000222.eps}} @>\div>>
  \raise-.15in\hbox{\includegraphics[width=.5in]{000210.eps}} @.\to 0
\end{CD}
\end{equation*}

We shall see many other discrete differential complexes below.

\section{Stability of Galerkin methods}

\vskip-5mm \hspace{5mm}

Consider first the solution of the Dirichlet problem for
Poisson's equation on a domain in $\R^n$:
\begin{equation*}
-\Delta u = f \text{ in $\Omega$},\quad
u=0 \text{ on $\partial\Omega$}.
\end{equation*}
The solution can be characterized as the minimizer of the energy
functional
\begin{equation*}
\mathcal E(u):=\frac12\int_\Omega|\grad u(x)|^2\,dx-\int_\Omega f(x)u(x)\,dx
\end{equation*}
over the Sobolev space $\mathaccent23 H^1(\Omega)$ (consisting of
$H^1(\Omega)$ functions vanishing on $\partial\Omega$), or as the solution
of the weak problem: find $u\in \mathaccent23 H^1(\Omega)$ such that
\begin{equation*}
\int_\Omega\grad u(x)\cdot\grad v(x)\,dx=\int_\Omega f(x)v(x)\,dx
\quad \text{for all $v\in\mathaccent23 H^1(\Omega)$}.
\end{equation*}
We may define an approximate solution $u_h$ by minimizing the Dirichlet
integral over a finite dimensional subspace $W_h$ of $\mathaccent23
H^1(\Omega)$; this is the classical Ritz method.  Equivalently, we may
use the Galerkin method, in which $u_h\in W_h$ is determined by the
equations
\begin{equation*}
\int_\Omega\grad u_h(x)\cdot\grad v(x)\,dx=\int_\Omega f(x)v(x)\,dx
\quad \text{for all $v\in W_h$}.
\end{equation*}
After choice of a basis in $W_h$ this leads to a system of linear
algebraic equations, and $u_h$ is computable.

Let $T_h$ denote the discrete solution operator $f\mapsto u_h$.  Then
it is easy to check that $T_h$ is bounded as a linear operator from
$H^{-1}(\Omega):=\mathaccent23H^1(\Omega)^*$ to $\mathaccent23H^1(\Omega)$
by a constant that depends only on the domain $\Omega$ (and, in
particular, doesn't increase if the space $W_h$ is enriched).
This says that the Galerkin method is \emph{stable}.  A consequence
is the \emph{quasioptimality estimate}
\begin{equation}\label{eq:quasiopt}
\|u-u_h\|_{H^1}\le c\inf_{v\in W_h}\|u-v\|_{H^1},
\end{equation}
for some constant $c$ depending only on the domain $\Omega$.
Note that there is no restriction on the subspace $W_h$ to obtain
this estimate.  Galerkin's method for a coercive elliptic problem
is always stable and convergence depends only on the approximation
properties of the subspace.  A natural choice for $W_h$ is the
Lagrange finite element space
of some degree $p$ with respect to some regular simplicial mesh of maximal
element size $h$, in which case Galerkin's method is a standard finite element
method.  In this case the right hand side of \eqref{eq:quasiopt}
is $O(h^p)$ provided that $u$ is sufficiently smooth.

Next consider the related eigenvalue problem, which arises
in the determination
of the fundamental frequencies of a drum.  That is, we seek
standing wave solutions $w(x,t)$ to the wave equation on some bounded domain
$\Omega\subset\R^2$ which vanish on $\partial\Omega$.
Assuming that the tension
and density of the drum membrane are unity, these solutions have the form
$w(x,t)=\alpha\cos(\sqrt{\lambda} t)u(x)+\beta\sin(\sqrt{\lambda} t)u(x)$
where $\alpha$ and $\beta$ are constants and $u$ and $\lambda$ satisfy
the eigenvalue problem
\begin{equation*}
-\Delta u = \lambda u \text{ in $\Omega$}, \quad u=0 \text{ on
$\partial\Omega$}.
\end{equation*}
The eigenvalues $\lambda$ form a sequence of positive numbers tending
to infinity.  The numbers $\sqrt\lambda/(2\pi)$ are the fundamental
frequencies of the drum and the functions $u$ give the corresponding
fundamental modes.

The eigenvalues and eigenfunctions are characterized variationally as
the critical values and critical points of the Rayleigh quotient
\begin{equation*}
\mathcal R(u)=\frac{\int_\Omega |\grad u(x)|^2\,dx}{\int_\Omega
|u(x)|^2\,dx},
\end{equation*}
defined for nonzero $u$ belonging to the Sobolev space
$\mathaccent23 H^1(\Omega)$.  The classical Rayleigh-Ritz method
for the approximation of eigenvalue problems determines
approximate eigenvalues $\lambda_h$ and eigenfunctions $u_h$ as
the  critical values and points of the restriction of $\mathcal R$
to the nonzero elements of some finite dimensional subspace $W_h$
of $\mathaccent23 H^1(\Omega)$. Equivalently, we can write the
eigenvalue problem in weak form: find $\lambda\in\R$ and nonzero
$u\in \mathaccent23 H^1(\Omega)$ such that
\begin{equation}\label{eq:eig}
\int_\Omega\grad u(x)\cdot\grad v(x)\,dx
=\lambda\int_\Omega u(x)v(x)\,dx
\quad \text{for all $v\in \mathaccent23 H^1(\Omega)$}.
\end{equation}
The Galerkin approximation of the eigenvalue problem, which is
equivalent to the Rayleigh-Ritz method, seeks $\lambda_h\in\R$ and
nonzero $u_h\in W_h$ such that
\begin{equation}\label{eq:eigh}
\int_\Omega\grad u_h(x)\cdot\grad v(x)\,dx
=\lambda_h\int_\Omega u_h(x)v(x)\,dx
\quad\text{for all $v\in W_h$}.
\end{equation}

We now discuss the convergence of this method.  Let $\lambda$ denote
the $j$th eigenvalue of the problem \eqref{eq:eig}. In the interest of
simplicity we assume that $\lambda$ is a simple eigenvalue, so the
corresponding eigenfunction $u$ is uniquely determined up to sign by
the normalization $\|u\|_{H^1}=1$.  Similarly let $\lambda_h$ and $u_h$
denote the $j$th eigenvalue of \eqref{eq:eigh}.  It can then be proved
(see, e.g., \cite{babuska-osborn} for much more general results)
that there exists a constant $c$ such that
\begin{equation}\label{eq:eigest}
\|u-u_h\|_{H^1}\le c\inf_{v\in W_h}\|u-v\|_{H^1}, \quad
|\lambda-\lambda_h|\le c\|u-u_h\|_{H^1}^2.
\end{equation}
In short, the eigenfunction approximation is quasioptimal
and the eigenvalue error is bounded by the square.
Again there is no restriction on the space $W_h$.

\begin{figure}[ht!]
\leftline{\hskip2.5in{\includegraphics[width=2.5in]{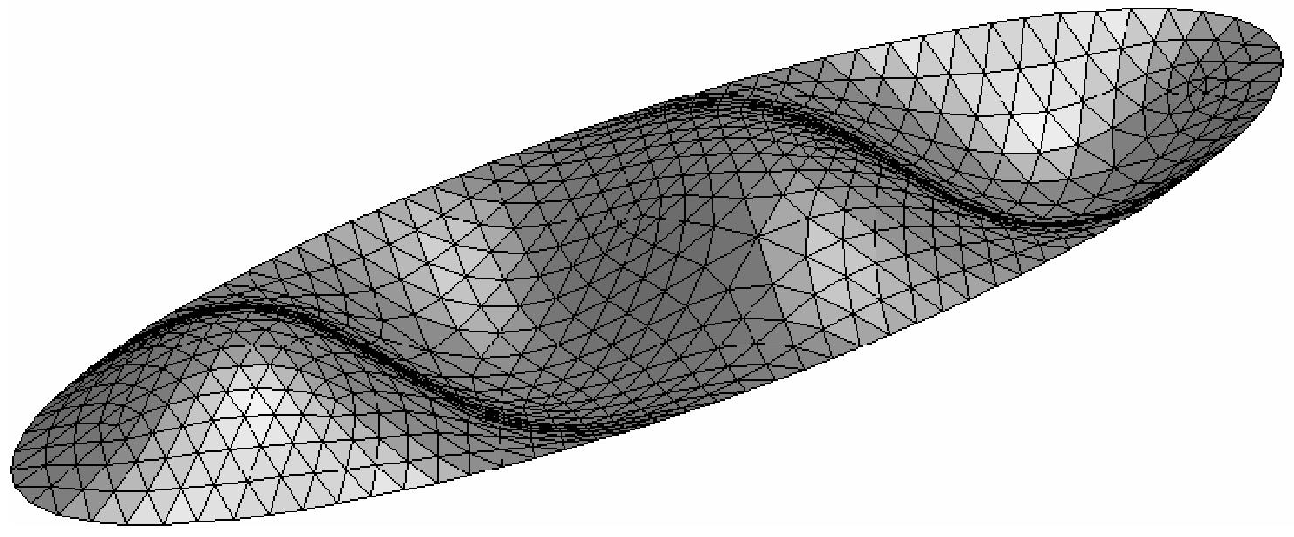}}} \vskip-.75in
\centerline{\raise.0in\hbox{\includegraphics[width=2.25in]{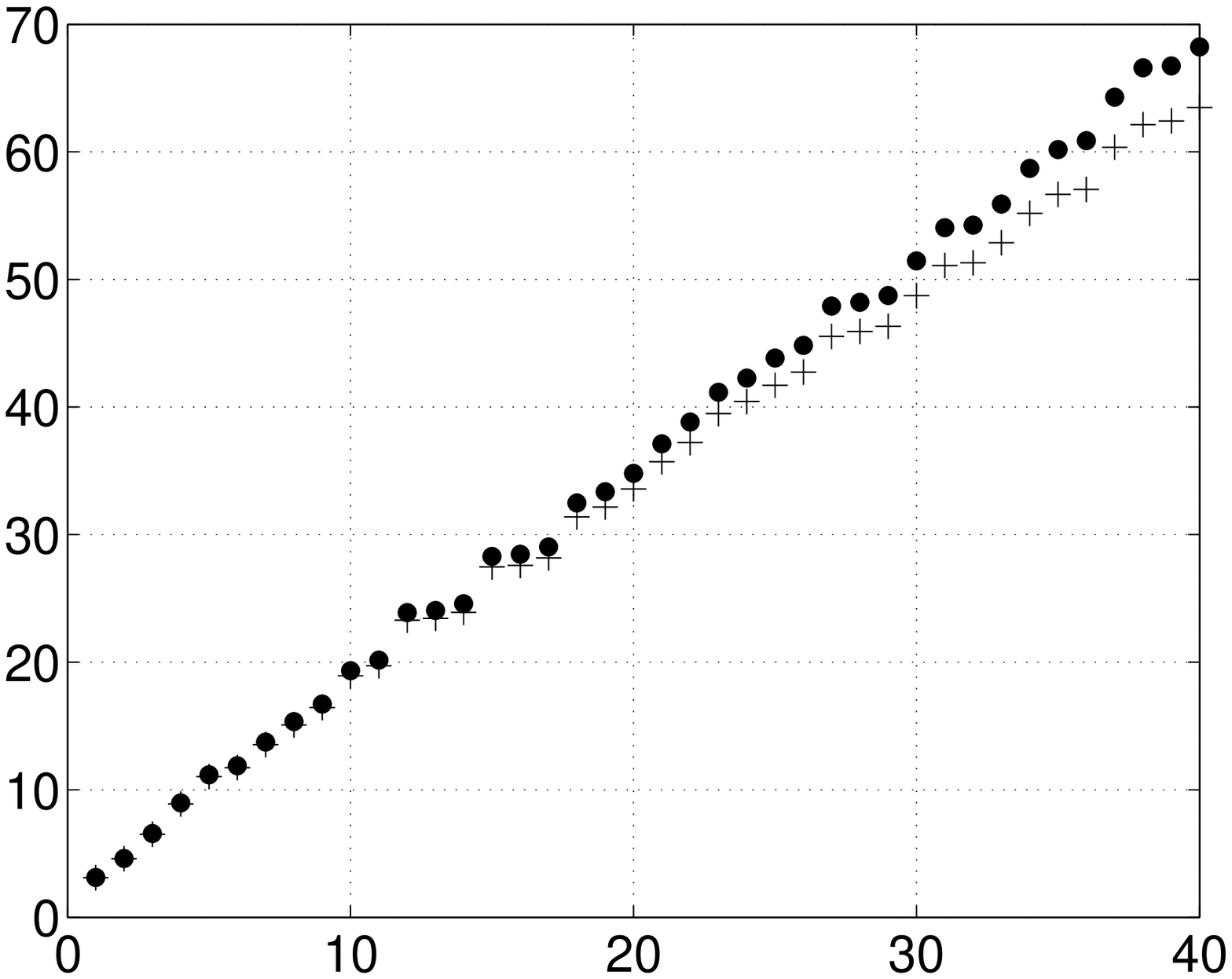}}\quad
\hbox{\includegraphics[width=2.5in]{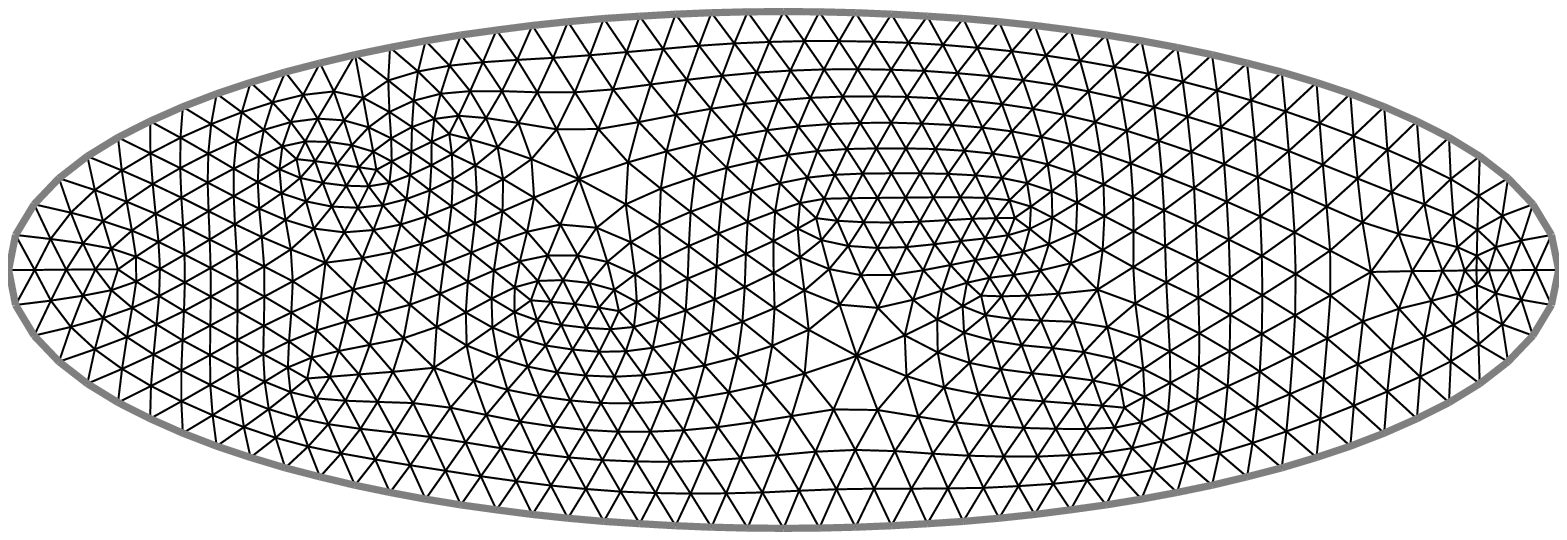}}}
\begin{center}
\begin{minipage}{10cm}
\caption{\label{fg:lapeig} The point plot shows the first 40
eigenvalues computed with piecewise linear finite elements with
respect to the triangulation shown ($\bullet$) versus the exact
eigenvalues ($+$).  The surface plot shows the computed
eigenfunction associated to the fourth eigenvalue. The mesh has
737 vertices, of which 641 are interior, and 1,376 triangles.}
\end{minipage}
\end{center}
\end{figure}

Figure~\ref{fg:lapeig} reports on the computation of
the eigenvalues of the Laplacian on an elliptical domain of aspect
ratio $3$ using Lagrange finite elements of degree $1$.

Now consider an analogous problem, the computation of the resonant
frequencies of an electromagnetic cavity occupying a region
$\Omega\subset\R^3$.  In this case we wish to find standing wave solutions
of Maxwell's equations.  If we take the
electric permittivity and the magnetic permeability to be unity and
assume a lossless cavity with perfectly conducting boundary, we are
led to the following eigenvalue problem for the  electric field:
find nonzero $E:\Omega\to\R^3$, $\lambda\in\R$ such that
\begin{equation}\label{eq:maxeig}
\curl\curl E=\lambda E, \quad \div E =0 \text{ in $\Omega$},
\quad E\times n =0 \text{ on $\partial\Omega$}.
\end{equation}
This is again an elliptic eigenvalue problem and the eigenvalues form
a sequence of positive numbers tending to infinity.  The divergence
constraint is nearly redundant in this eigenvalue problem.  Indeed
if $\curl\curl E=\lambda E$ for $\lambda>0$, then $\div
E=\lambda^{-1}\div\curl\curl E =0$ since the divergence of a curl
vanishes.  Thus the eigenvalue problem
\begin{equation}\label{eq:maxeig0}
\curl\curl E=\lambda E \text{ in $\Omega$},
\quad E\times n =0 \text{ on $\partial\Omega$},
\end{equation}
has the same eigenvalues and eigenfunctions as \eqref{eq:maxeig} except
that it also admits $\lambda=0$ as an eigenvalue, and the corresponding
eigenspace is infinite-dimensional (it contains the gradients of all
smooth functions vanishing on the boundary of $\Omega$).  The
eigenvalues and eigenfunctions are now critical points and values of the
Rayleigh quotient
\begin{equation*}
\mathcal R(E)=\frac{\int_\Omega |\curl E(x)|^2dx}{\int_\Omega|E(x)|^2dx},
\end{equation*}
over the space of nonzero fields $E$ in $\mathaccent 23 H(\curl,\Omega)$, which is
defined to be the space of functions for which both the above
integrals exist and are finite and which have vanishing tangential
component on the boundary (i.e., $E\times n=0$ on $\partial\Omega$).

\begin{figure}[ht!]
\centerline{\includegraphics[height=2in]{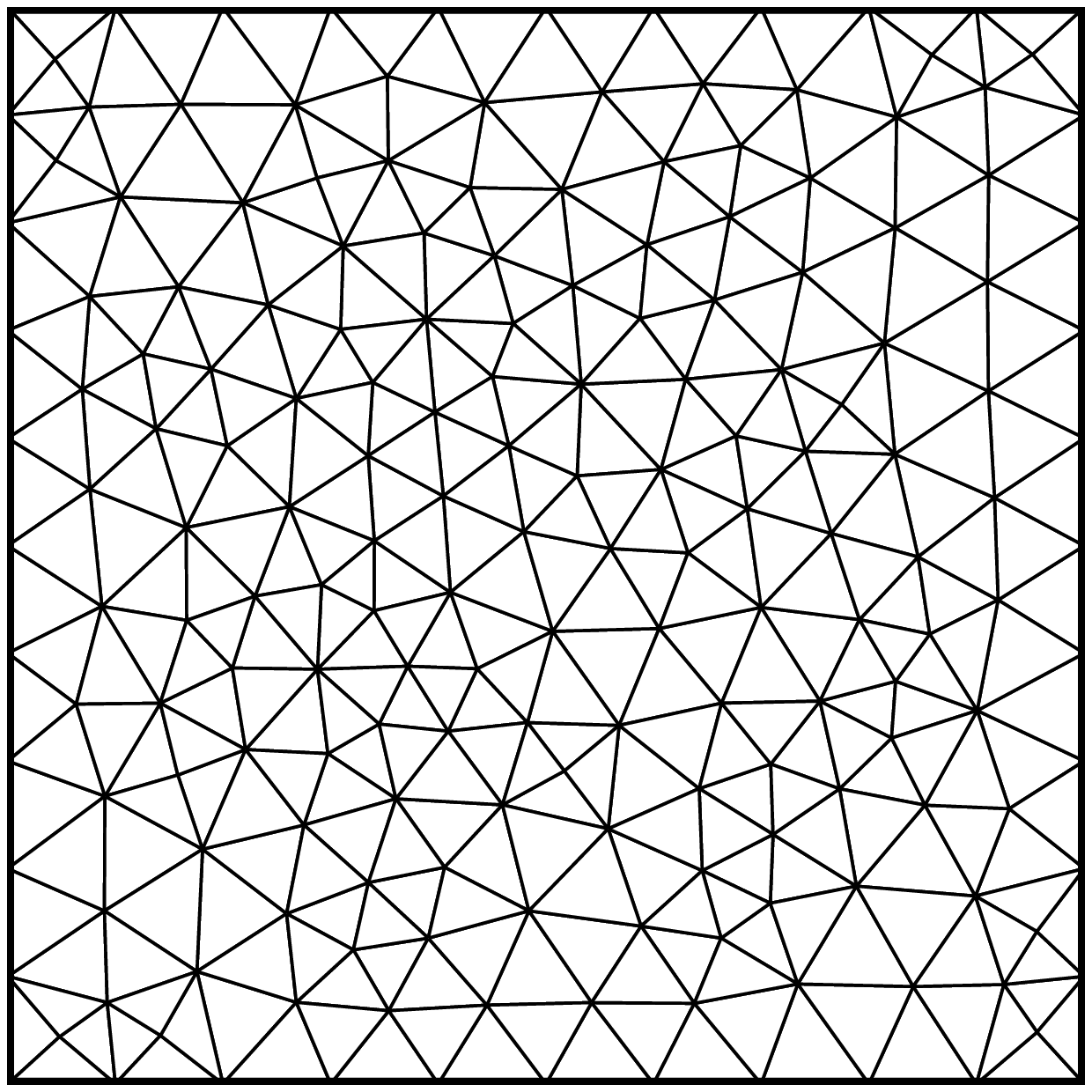}\qquad
\includegraphics[height=2in]{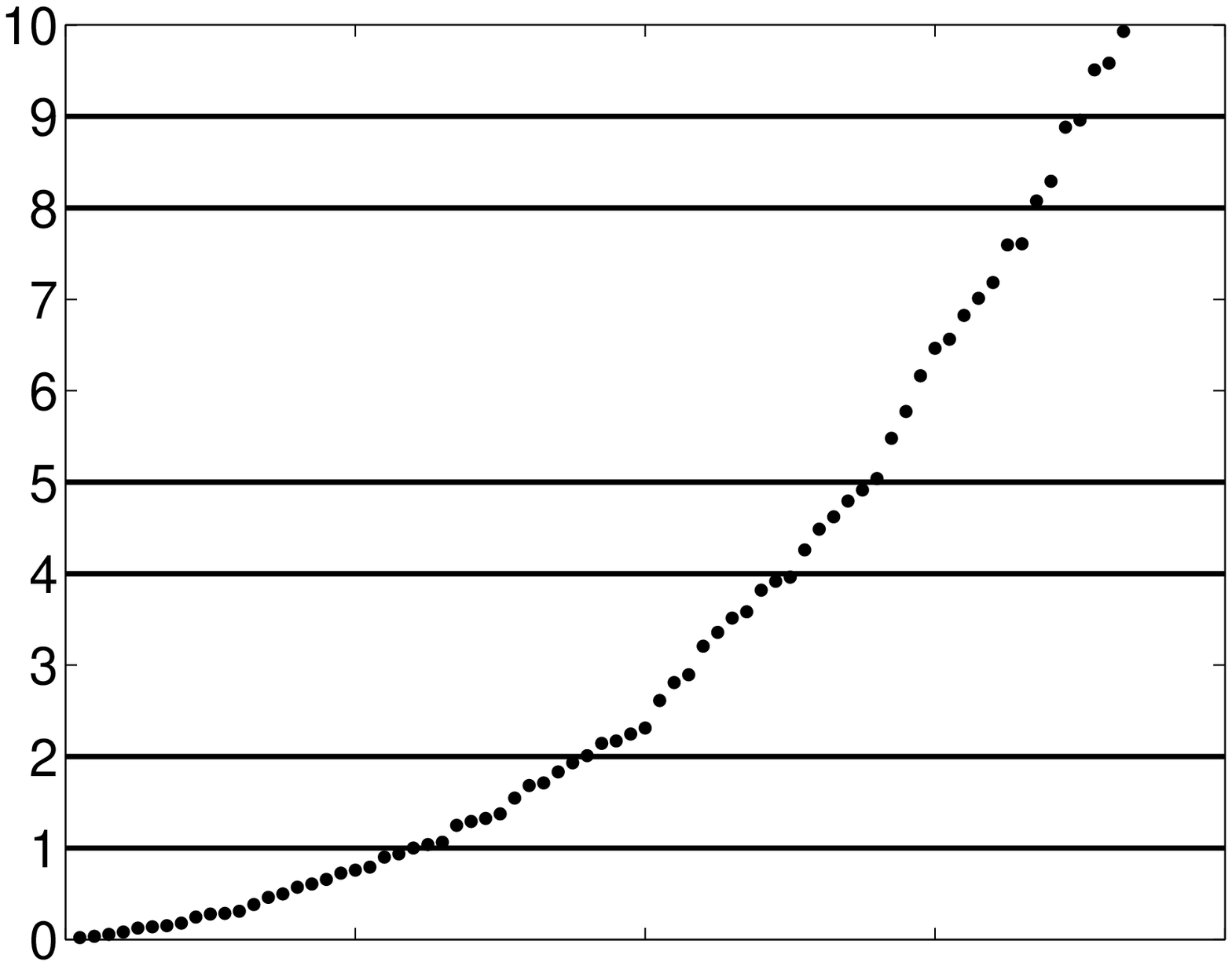}}
\begin{center}
\begin{minipage}{10cm}
\caption{\label{fg:nodeig} The plot shows the first 73 eigenvalues
computed with piecewise linear finite elements for the resonant
cavity problem on the square using the mesh shown.  They bear no
relation to the exact eigenvalues, $1$, $1$, $2$, $4$, $4$, \dots,
indicated by the horizontal lines.}
\end{minipage}
\end{center}
\end{figure}

In Figure~\ref{fg:nodeig} we show the result of approximating a
two-dimensional version of this eigenvalue problem using the
Rayleigh-Ritz method or, equivalently, the Galerkin method with
continuous piecewise linear vector fields on $\Omega$ whose
tangential components vanish on the boundary (the first element
depicted in Figure~\ref{fg:eltsvec}).  For $\Omega$ we take a
square of side length $\pi$, in which case the nonzero eigenvalues
are known to be all numbers of the form $\lambda=m^2+n^2$ with
$0\le m,n\in \mathbb Z$ not both zero, and the corresponding
eigenfunctions are $E=(\sin m y,\sin n x)$.   For the mesh
pictured, the finite element space has dimension $290$.  We find
that $73$ of the $290$ computed eigenvalues are between $0$ and
$10$ and that they have no tendency to cluster near the integers
$1,1,2,4,4,5,5,8,9,9$ which are the exact eigenvalues between $0$
and $10$.  Thus this numerical method is useless: the computed
eigenvalues bear no relation to the true eigenvalues!  The
analogue of \eqref{eq:eigest} is surely not true.

\begin{figure}[ht!]
\centerline{\includegraphics[height=2in]{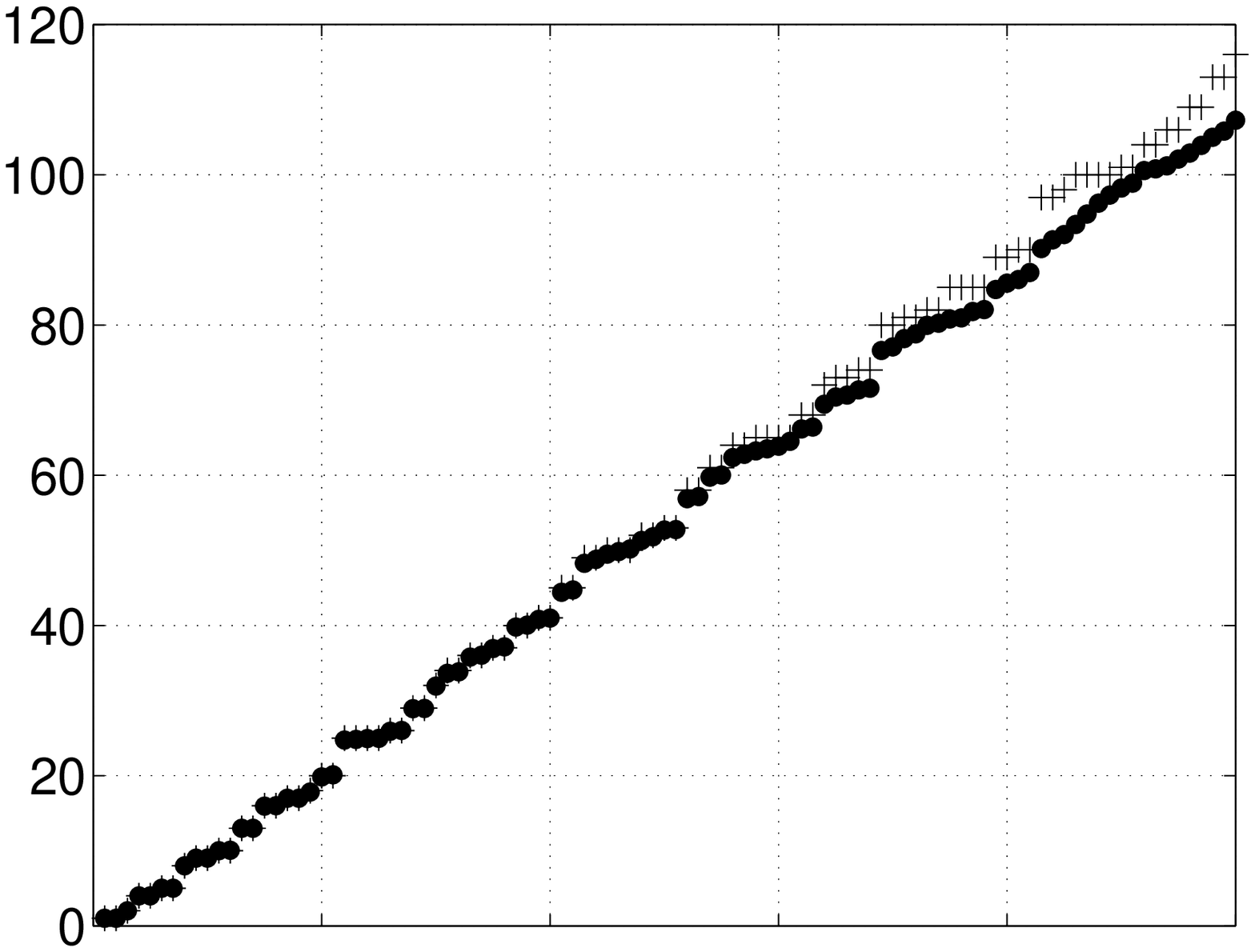}\qquad
\raise.05in\hbox{\includegraphics[height=1.9in]{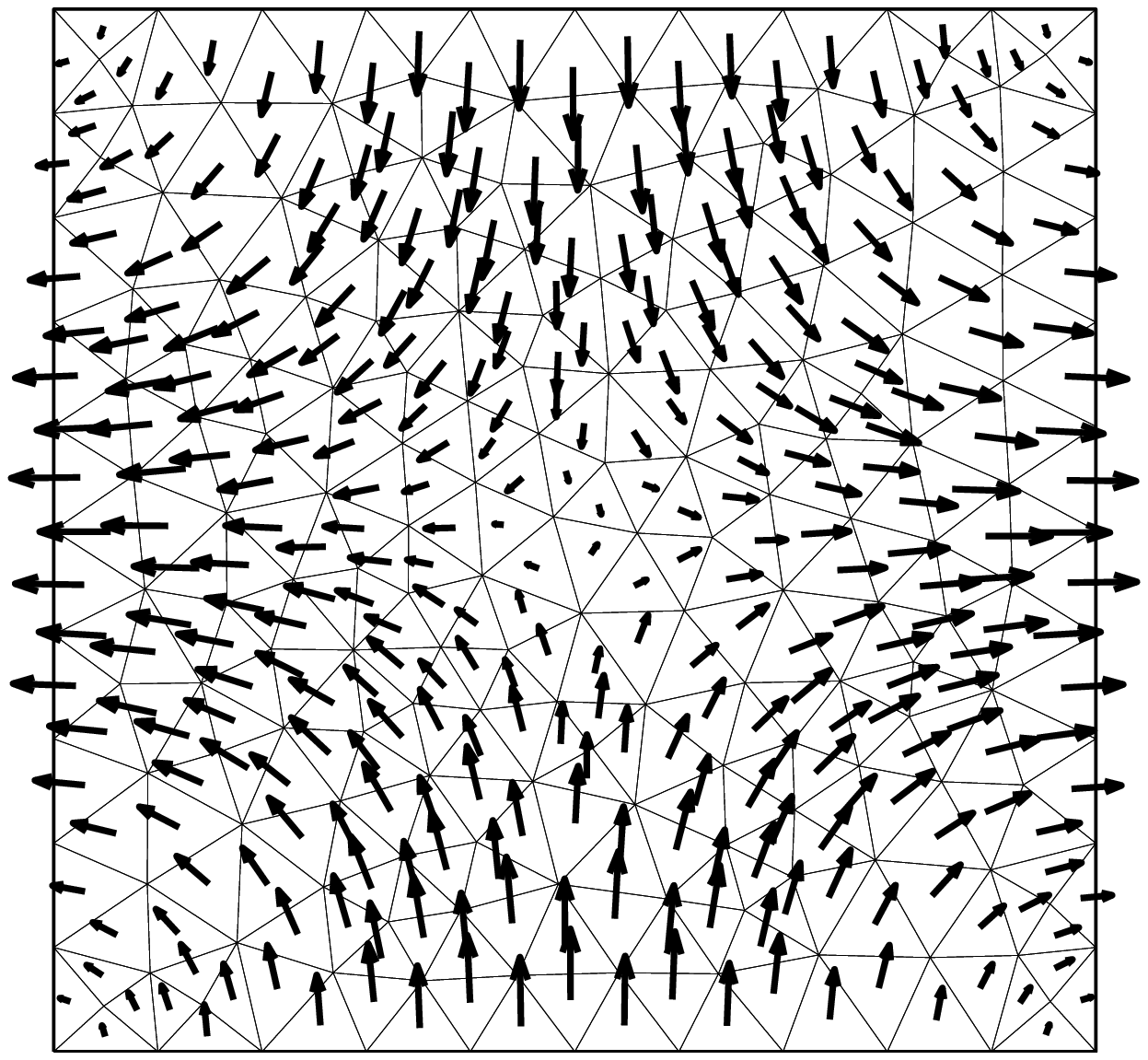}}} \vskip-1.9in
\leftline{\hskip.25in\includegraphics[height=.9in]{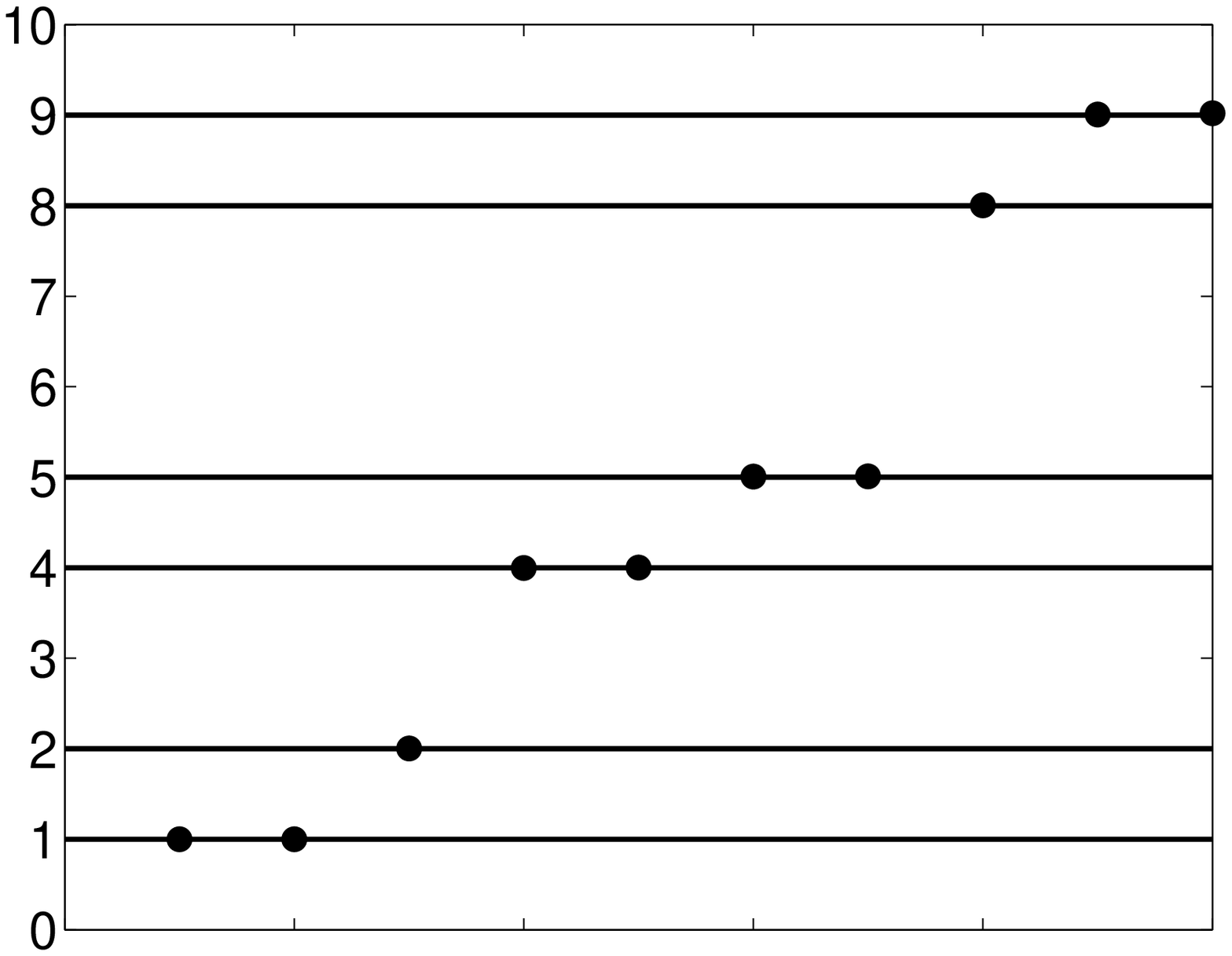}} \vskip.9in
\begin{center}
\begin{minipage}{10cm}
\caption{\label{fg:edgeeig} The first plot shows the first 100
positive eigenvalues for the resonant cavity problem on the square
computed with lowest order edge elements using the mesh of
Figure~\ref{fg:nodeig}.  The error in the first $54$ eigenvalues
is below $2\%$.    The inset focuses on the first $10$
eigenvalues, for which the error is less than $0.25\%$.  The
second plot shows the vector field associated to the third
positive eigenvalue.}
\end{minipage}
\end{center}
\end{figure}

If instead we choose the lowest order edge elements as the finite
element space (Figure~\ref{fg:eltsvec}, top right), we get very
different results.  Using the same mesh, the edge finite element space has
dimension $472$.  It turns out that $145$ of the computed eigenvalues are
zero (to within round-off), and
the subsequent eigenvalues are $0.9998$, $0.9999$, $2.0023$, $3.9968$,
$4.0013$, \dots, i.e., excellent approximations of the exact
eigenvalues.  See Figure~\ref{fg:edgeeig}.

The striking difference between the behavior of the continuous piecewise
linear finite elements and the edge elements for the resonant cavity
problem is a question of stability.  We shall return to this below,
after examining stability in a simpler context.

\section{Stability of mixed formulations}

\vskip-5mm \hspace{5mm}

Consider now the Dirichlet problem
\begin{equation*}
-\div C\grad u =f\text{ in $\Omega$}, \quad u=0 \text{ on
$\partial\Omega$},
\end{equation*}
where $\Omega$ is a domain in $\R^3$ and the coefficient $C$ is a
symmetric positive definite matrix at each point.
We may again characterize $u$ as a minimizer of the energy functional
\begin{equation*}
u\mapsto \frac 12\int C\grad u\cdot\grad u dx-\int fu\,dx
\end{equation*}
and use the Ritz method.  This procedure is always stable.

However, for some purposes it is preferable to work with the equivalent
first order system
\begin{equation}\label{eq:fo}
\sigma=C\grad u, \quad -\div\sigma=f.
\end{equation}
The pair $(\sigma,u)$ is then characterized variationally as the
unique critical point of the functional
\begin{equation}\label{eq:mixedfun}
\mathcal L(\sigma,u)=
\int_\Omega(\frac12 C^{-1}\sigma\cdot\sigma+u\div\sigma)dx-\int_\Omega  fu\,dx
\end{equation}
over $H(\div,\Omega)\times L^2(\Omega)$.
Note that $(\sigma,u)$ is a saddle-point of $\mathcal L$,
not an extremum.  Numerical discretizations based on such saddle-point
variational principles are called \emph{mixed methods}.

It is worth interpreting the system \eqref{eq:fo} in the language of
differential forms, because this brings some insight.  The function $u$
is a $0$-form, and the operation $u\mapsto\grad u$ is just exterior
differentiation.  The vector field $\sigma$ is a proxy for a $2$-form
and the operation $\sigma\mapsto\div\sigma$ is again exterior
differentiation.  The loading function $f$ is the proxy for a
$3$-form.  Since $\grad u$ is the proxy for a $1$-form, it must be that
the operation on differential forms that
corresponds to multiplication by $C$ takes $1$-forms to $2$-forms.
In fact, if we untangle the identifications, we find that
multiplication by $C$ is a Hodge star operation.  A Hodge star operator
defines an isomorphism of $\bigwedge^k(\Omega)$ onto
$\bigwedge^{3-k}(\Omega)$. To determine a particular such operator, we
must define an inner product on the tangent space $\R^3$ at each
point of $\Omega$.  The positive
definite matrix $C$ does exactly that. Many of the partial differential
equations of mathematical physics admit similar interpretations in
terms of differential forms. For a discussion of this in the context of
discretization, see \cite{hiptmair}.

A natural approach to discretization of the mixed variational principle
is to choose subspaces $S_h\subset H(\div,\Omega)$, $V_h\subset
L^2(\Omega)$ and seek a critical point $(\sigma_h,u_h)\in S_h\times
V_h$.  This is of course equivalent to a Galerkin method and leads to a
system of linear algebraic equations. However in this case,
\emph{stability is not automatic}. It can happen that the discrete
system is singular, or more commonly, that the norm of the discrete
solution operator grows unboundedly as the mesh is refined.

In a fundamental paper, Brezzi \cite{brezzi} established two
conditions that together are sufficient (and essentially necessary) for
stability.  Brezzi's theorem applied to a wide class of saddle-point
problems, but for simplicity we will state the stability conditions for the
saddle-point problem associated to the functional \eqref{eq:mixedfun}.
\begin{itemize}
\item[(S1)] There exists $\gamma_1>0$ such that
\begin{equation*}
\int_\Omega C^{-1}\tau\cdot\tau\,dx\ge\gamma_1\|\tau\|_{H(\div)}^2,
\end{equation*}
for all $\tau\in S_h$ such that $\int\div\tau\, v\,dx=0$  for all
$v\in V_h$.
\item[(S2)] There exists $\gamma_2>0$ such that for all $v\in V_h$
there exists nonzero $\tau\in S_h$ satisfying
\begin{equation*}
\int_\Omega v\div\tau\,dx \ge \gamma_2\|v\|_{L^2}\|\tau\|_{H(\div)}.
\end{equation*}
\end{itemize}

\noindent
{\bf Theorem (Brezzi)} {\it If the stability conditions {\rm (S1)} and
{\rm (S2)} are satisfied,
then $\mathcal L$ admits a unique critical point $(\sigma_h,u_h)$ over $S_h\times
V_h$, the solution operator $f\mapsto (\sigma_h,u_h)$ is bounded
$L^2(\Omega)\to H(\div,\Omega)\times L^2(\Omega)$, and the quasioptimal estimate
\begin{equation*}
\|\sigma-\sigma_h\|_{H(\div)}+\|u-u_h\|_{L^2}\le
c\inf_{(\tau,v)\in S_h\times
V_h}(\|\sigma-\tau\|_{H(\div)}+\|u-v\|_{L^2})
\end{equation*}
holds with $c$ depending on $\gamma_1$ and $\gamma_2$.}

The stability conditions of Brezzi strongly limit the choice of the
mixed finite element spaces $S_h$ and $V_h$. Condition (S1) is
satisfied if the indicated functions $\tau\in S_h$, those whose
divergence is orthogonal to $V_h$, are in fact divergence-free.  (In
practice, this is nearly the only way it is satisfied.) This certainly
holds if $\div S_h\subset V_h$, and so such as inclusion is a common
design principle of mixed finite element spaces. On the other hand,
condition (S2) is most easily satisfied if $\div S_h\supset V_h$,
because in this case, given $v\in V_h$, we can choose $\tau\in S_h$
with $\div\tau=v$, so $\int_\Omega v\div\tau\,dx=\|v\|_{L^2}^2$, and
the second condition will be satisfied as long as we can insure that
$\|\tau\|_{H(\div)}\le\gamma_2^{-1}\|v\|_{L^2}$.  In short, we need to
know that $\div$ maps $S_h$ onto $V_h$ and that $\div|_{S_h}$ admits a
bounded one-sided inverse.

The face elements of Raviart-Thomas and Nedelec were designed to
satisfy both these conditions.  Specifically, let $S_h$ again
denote the space of face elements of lowest degree (whose element
diagram is shown in the middle of the second row of
Figure~\ref{fg:eltsvec}), and $V_h$ the space of piecewise
constants.\footnote{It may seem odd to seek $u_h$ in $V_h$, a
space of discrete $3$-forms, rather than in a space of $0$-forms,
since $u$ is a $0$-form. The resolution is through a Hodge star
operator, this time formed with respect to the Euclidean inner
product on $\R^3$.  In the mixed method $u_h$ is a discrete
$3$-form, approximating the image of $u$ under this star
operator.}  We know that $S_h\subset H(\div,\Omega)$ so these
elements are admissable for the mixed variational principle.
Moreover, we have $\div S_h\subset V_h$, so (S1) holds.

To verify (S2), we refer to the commutative diagram
\eqref{eq:cd}. Given $v\in V_h$, we can solve the Poisson equation
$\Delta\phi=v$ and take $\sigma=\grad\phi$ to obtain a function with
$\div\sigma=v$ and $\|\sigma\|_{H^1}\le C\|v\|_{L^2}$.  Now let
$\tau=\Pi^S_h\sigma\in S_h$.  Then
\begin{equation*}
\div\tau=\div\Pi^S_h\sigma=\Pi^V_h\div\sigma
=\Pi^V_hv =v,
\end{equation*}
where we have used the commutativity and the fact that
$v\in V_h$.  Moreover $\|\tau\|_{H(\div)}\le c\|\sigma\|_{H^1}\le
c'\|v\|_{L^2}$, where we used the boundedness
of $\Pi^S_h$ on $H^1(\Omega,\R^3)$.  This shows that $\div V_h=S_h$ and
establishes a bound on the one-sided inverse, and
so verifies (S2).  Of course, the same
argument shows the stability of a mixed method
based on higher order face elements as well.

Thus we see that the stability of the mixed finite element method
depends on the properties of the spaces $V_h$ and $S_h$ encoded
in the rightmost square of the commutative diagram  \eqref{eq:cd}.

Now let us return to the resonant cavity eigenvalue problem
\eqref{eq:maxeig0} for which we explored the Galerkin method: find $\lambda_h\in\R$,
$0\ne E_h\in Q_h$ such that
\begin{equation}\label{eq:gal}
\int_\Omega \curl E_h\cdot\curl F\,dx = \lambda_h\int_\Omega E_h\cdot
F\,dx\quad \text{for all $F\in Q_h$}.
\end{equation}
We saw that  if $Q_h\subset\mathaccent23H(\curl,\Omega)$ is taken to be a space of edge elements this
method gives good results  in that the positive eigenvalues of the
discrete problem are good approximations for the positive eigenvalues
of the continuous problem.  However, the simple choice of  Lagrange
finite elements did not give good results.  We now explain the good
performance of the edge elements based on the middle square of the
commutative diagram \eqref{eq:cd}.  Following Boffi et.~al \cite{boffi}
we set $P_h=\curl Q_h$ and
introduce the following mixed discrete eigenvalue problem: find
$\lambda_h\in\R$, $0\ne (E_h,p_h)\in Q_h\times P_h$ such that
\begin{gather}\label{eq:mixed}
\int_\Omega E_h\cdot F\,dx + \int_\Omega \curl F\cdot p_h\,dx = 0
\quad\text{for all $F\in Q_h$},\\
\int_\Omega \curl E_h\cdot q\,dx = -\lambda_h\int_\Omega p_h\cdot
q\,dx\quad\text{for all $q\in P_h$}.
\end{gather}
It is then easy to verify that if $\lambda_h$, $E_h$ is a solution
to \eqref{eq:gal} with $\lambda_h>0$, then $\lambda_h$,
$(E_h,\lambda_h^{-1}\curl E_h)$ is a solution to \eqref{eq:mixed}, and if
$\lambda_h$, $(E_h,p_h)$ is a solution to \eqref{eq:mixed} then
$\lambda_h>0$ and $\lambda_h$, $E_h$ is a solution to \eqref{eq:gal}.
In short, the two problems are equivalent except that the former
admits a zero eigenspace which the mixed
formulation suppresses.  As explained in
\cite{boffi}, the accuracy of the mixed eigenvalue problem
\eqref{eq:mixed} hinges on the stability of the corresponding mixed
source problem.  This is a saddle-point problem of the sort studied
by Brezzi,  and so stability depends on conditions analogous to (S1) and
(S2).  The proof of these
conditions in case $Q_h$ is the space of edge elements follows,
as in the preceding stability verification, from surjectivity and
commutativity properties encoded in the diagram \eqref{eq:cd}.

The diagram can also be used to explain the zero eigenspace computed
with edge elements.  Recall that in the case of the mesh shown in
Figure~\ref{fg:nodeig}, this space had dimension $145$.  In fact, this
eigenspace is simply the null space of the curl operator restricted
to $Q_h$.  Referring again to the commutative diagram \eqref{eq:cd},
this is the gradient of the space $W_h$ of linear Lagrange elements
vanishing on the boundary.  Its dimension is therefore exactly the
number of interior nodes of the mesh.

\section{The elasticity complex}

\vskip-5mm \hspace{5mm}

Let $\mathbb S$ denote the space of $3\times3$ symmetric matrices.
Given a volumetric loading density $f:\Omega\to\R^3$, the system of linearized
elasticity determines the displacement field $u:\Omega\to\R^3$ and the stress
field $\sigma:\Omega\to\mathbb S$ induced in the elastic domain $\Omega$
by the equations
\begin{equation*}
\sigma= C\eps u, \quad -\div\sigma=f,
\end{equation*}
together with boundary conditions such as $u=0$ on $\partial\Omega$.
Here $\eps u$ is the symmetric part of the matrix $\grad u$, and
the elasticity tensor $C:\mathbb S\to\mathbb S$ is a symmetric positive
definite linear operator describing the particular elastic material,
possibly varying from point to point.

The solution $(\sigma,u)$ may be characterized variationally as a
saddle-point of the Hellinger-Reissner functional
\begin{equation}\label{eq:mixedelas}
\mathcal L(\sigma,u)=\int_\Omega(\frac12 C^{-1}\sigma:\sigma +
u\cdot\div\sigma)dx - \int_\Omega f\cdot u\,dx
\end{equation}
over $H(\div,\Omega,\mathbb S)\times L^2(\Omega,\R^2)$ (i.e.,
$\sigma$ is sought in the space of square-integrable
symmetric-matrix-valued functions whose divergence by rows
is square-integrable, and $u$ is sought among all square-integrable
vector fields).

For a mixed finite element method, we need to specify finite element
subspaces $S_h\subset H(\div,\Omega,\mathbb S)$ and $V_h\subset
L^2(\Omega,\R^2)$ and restrict the domain of the variational problem.
Of course the spaces must be carefully designed if the mixed method
is to be stable: the analogues of the stability conditions (S1) and (S2)
must be satisfied.  The functional \eqref{eq:mixedelas} is quite similar
in appearance to \eqref{eq:mixedfun} and so it might be expected that
the mixed finite elements developed for the latter (the face elements
for $\sigma$ and discontinuous elements for $u$) could be adapted to
the case of elasticity.  In fact, the requirement of symmetry of the
stress tensor and, correspondingly, the replacement of the gradient by the
symmetric gradient, changes the structure significantly.  Four decades
of searching for mixed finite elements for elasticity beginning
in the  1960s did not yield any stable elements with polynomial shape
functions.

Using discrete differential complexes, R.~Winther and the author
recently developed the first such elements for elasticity problems in
two dimensions \cite{mixedelas}.  (The three-dimensional case remains open.)
For elasticity, the
displacement and stress fields cannot be naturally interpreted as
differential forms and the relevant differential complex is not the de~Rham complex.
In three dimensions it is instead the \emph{elasticity complex}:
\begin{equation*}
\begin{CD}
\mathbb T\hookrightarrow\,@. C^\infty(\Omega,\R^3) @>\epsilon>>
C^\infty(\Omega,\mathbb S) @>J>> C^\infty(\Omega,\mathbb S) @>\div>>
C^\infty(\Omega,\R^3) @.\to0.
\end{CD}
\end{equation*}
Here the operator $J$ is a \emph{second order differential
operator} which acts on a symmetric matrix field by first
replacing each row with its curl and then replacing each column
with its curl to obtain another symmetric matrix field.  The
resolved space $\mathbb T$ is the six-dimensional space of
infinitesimal rigid motions, i.e., the same space of linear
polynomials $a+b\times x$ which arose as the shape functions for
the lowest order edge elements.  If the domain $\Omega$ is
topologically trivial, this complex is exact.  Although it
involves a second order differential operator, and so looks quite
different from the de~Rham complex, Eastwood \cite{eastwood}
recently pointed out that it can be derived from the de~Rham
complex via a general construction known as the
Bernstein-Gelfand-Gelfand resolution.

In two dimensions the elasticity complex takes the form
\begin{equation*}
\begin{CD}
\mathbb P_1\hookrightarrow\,@.
C^\infty(\Omega) @>J>> C^\infty(\Omega,\mathbb S) @>\div>>
C^\infty(\Omega,\R^2) @.\to0,
\end{CD}
\end{equation*}
where now the second order differential operator is
\begin{equation*}
J=\begin{pmatrix} \displaystyle\frac{\partial^2}{\partial x_2^2} &
 -\displaystyle\frac{\partial^2}{\partial x_1\partial x_2}\\\\
 -\displaystyle\frac{\partial^2}{\partial x_1\partial x_2} &
 \displaystyle\frac{\partial^2}{\partial
 x_1^2}\end{pmatrix}.
\end{equation*}

\begin{figure}[ht!]
\centerline{\raise-.17in\hbox{\includegraphics[width=1in]{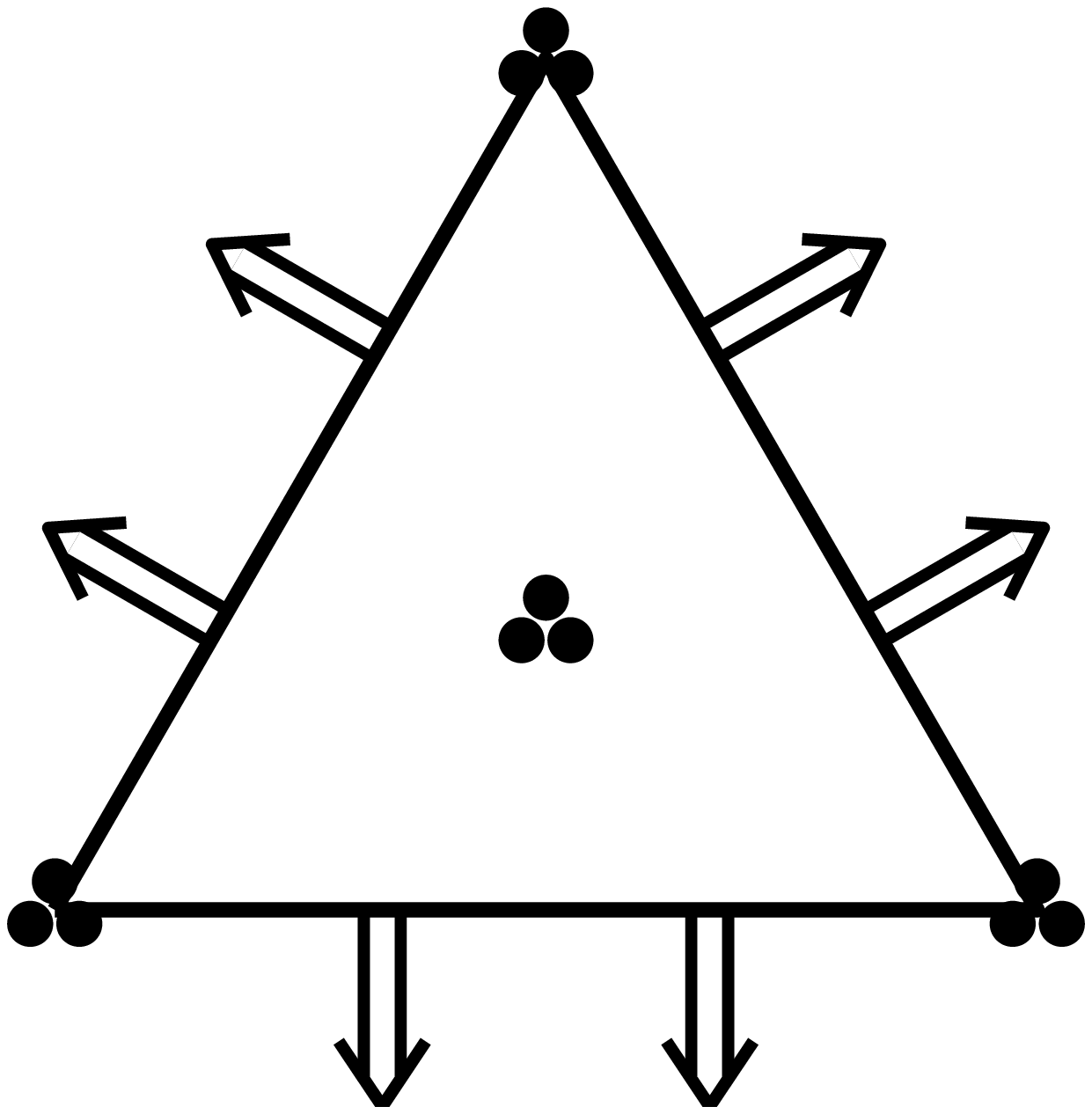}}\quad
\includegraphics[width=1in]{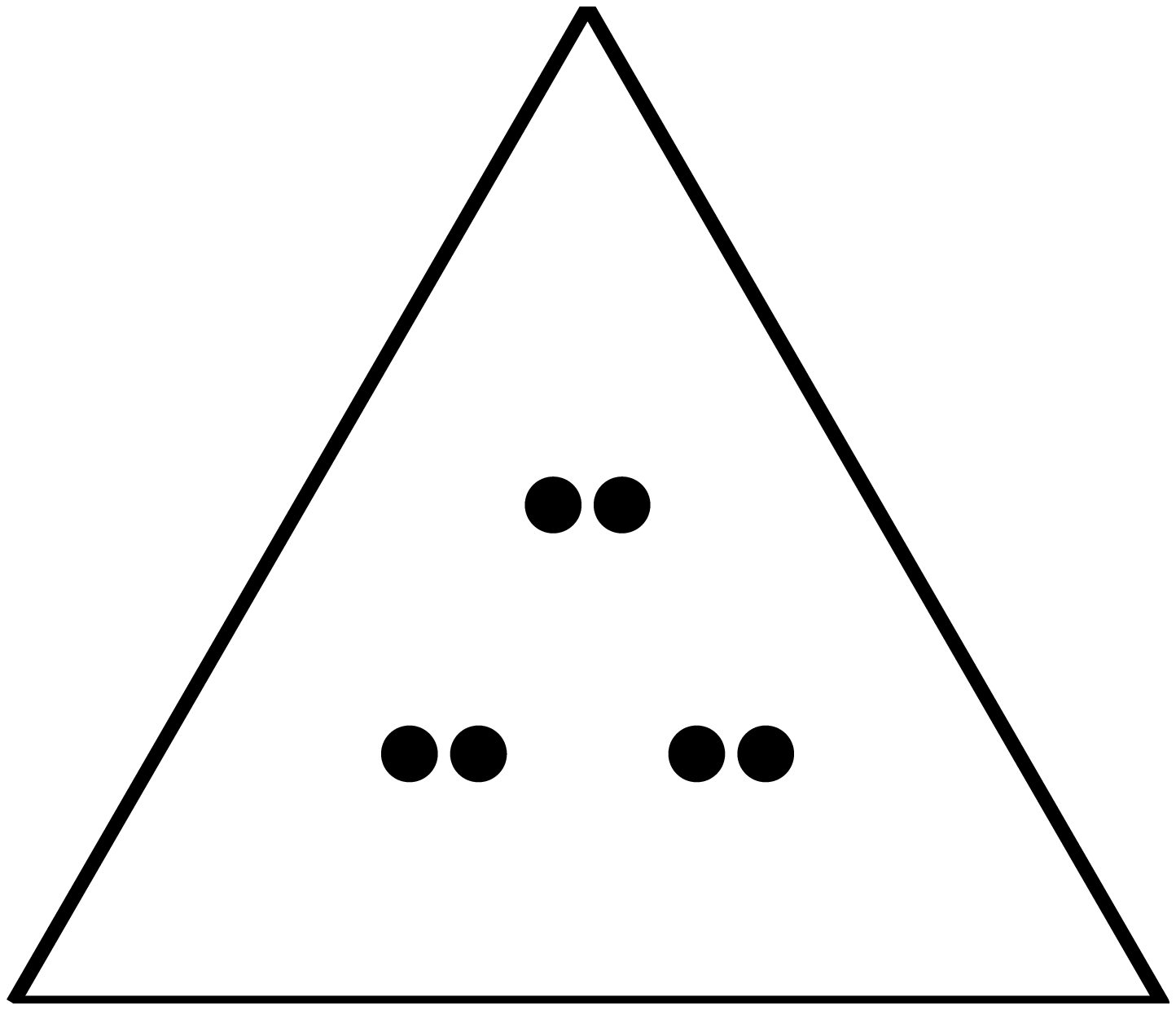}}
\begin{center}
\begin{minipage}{10cm}
\caption{\label{fg:aw} Element diagram for the new mixed finite
elements for elasticity, lowest order case.}
\end{minipage}
\end{center}
\end{figure}
In the lowest order case, the finite elements we introduced in
\cite{mixedelas}, for which the element diagrams can be
seen in Figure~\ref{fg:aw}, use discontinuous piecewise linear vector fields
for the displacement field and a piecewise polynomial space which we shall
now describe for the stress field.  The shape functions on an arbitrary
triangle $T$ are given by
\begin{equation*}
S_T = \{\,\tau\in \mathbb P_3(T,\mathbb S)\,|\, \div\tau\in \mathbb
P_1(T,\R^2)\,\},
\end{equation*}
which is a $24$-dimensional space consisting of all quadratic symmetric
matrix fields on $T$ together with the divergence-free cubic fields.
The degrees of freedom are
\begin{itemize}
\item the values of three components of $\tau(x)$ at each vertex $x$ of $T$
(9 degrees of freedom)
\item the values of the moments of degree $0$
and $1$ of the two components of $\tau n$ on each edge $e$ of $T$
(12 degrees of freedom)
\item the value of the three components of the moment of degree $0$
of $\tau$ on $T$ (3 degrees of freedom)
\end{itemize}
Note that these degrees of freedom are enough to ensure continuity
of $\tau n$ across element faces, and so will furnish a finite
element subspace of $H(\div,\Omega,\mathbb S)$.  The continuity
is not however, the minimal needed for inclusion in $H(\div)$.
The degrees of freedom also enforce continuity at the vertices,
which is not required for membership in $H(\div)$.  For various
reasons, it would be
useful to have a mixed finite element for elasticity that does
not use vertex degrees of freedom.  But, as we remark below,
this is not possible if we restrict to polynomial shape functions.

In order to have a well-defined finite element, we must verify that the
$24$ degrees of freedom form a basis for the dual space
of $S_T$.  We include this verification since it illustrates an aspect of the
role of the elasticity complex.
Since $\dim S_T=24$, we need only show that if all the
degrees of freedom vanish for some $\tau\in S_T$, then $\tau=0$.
Now $\tau n$ varies cubically along each edge, vanishes at the
endpoints, and has vanishing moments of degree $0$ and $1$.  Therefore
$\tau n\equiv 0$. Letting $v=\div\tau$, a linear vector field on $T$, we get
by integration by parts that
\begin{equation*}
\int_T v^2\,dx = -\int_T\tau:\eps v\,dx + \int_{\partial T} \tau n\cdot v\,ds=0
\end{equation*}
since the integral of $\tau$ vanishes as well as $\tau n$.
Thus $\tau$ is divergence-free.  In view of the exactness of the
elasticity complex,
$\tau=Jq$ for some smooth function $q$.  Since all the second partial
derivatives of $q$ belong to $\mathbb P_3(T)$, $q\in\mathbb P_5(T)$.
Adjusting by an element of $\mathbb P_1(T)$ (the null space of $J$), we may
take $q$ to vanish at the vertices.  Now $\partial^2q/\partial s^2 = \tau n\cdot n=0$
on each edge, whence $q$ is identically zero on $\partial T$.  This implies
that the gradient of $q$ vanishes at the vertices.  Since
$\partial^2 q/\partial s\partial n = -\tau n \cdot t=0$ on each edge (with $t$ a unit vector
tangent to the edge), we conclude that $\partial q/\partial n$ vanishes identically
on $\partial T$ as well.  Since $q$ has degree at most $5$, it must vanish
identically.

Let $\Pi^S_h: C^\infty(\Omega,\mathbb S)\to S_h$ denote the projection
associated with the supplied degrees of freedom, and
$\Pi^V_h:C^\infty(\Omega,\R^2)\to V_h$ the $L^2$-projection.  For
any triangle $T$, $\tau\in C^\infty(\Omega,\mathbb S)$, and $v\in \mathbb
P_1(\Omega,\R^2)$, we have
\begin{equation*}
\int_T\div(\tau-\Pi^S_h\tau)\cdot v\,dx
=-\int_T(\tau-\Pi^S_h\tau):\eps v\,dx +\int_{\partial
T}(\tau-\Pi^S_h\tau)n\cdot v\,ds.
\end{equation*}
The degrees of freedom entering the definition of $\Pi^S_h$ ensure that
the right hand side vanishes, and from this we obtain the commutativity
$\div\Pi^S_h\tau=\Pi^V_h\div\tau$ which is essential for stability.
(Actually a technical difficulty arises here, since $\Pi^S_h$ as given is
not bounded on $H^1(\Omega,\mathbb S)$.  See \cite{mixedelas} for the
resolution.)  Note that, by their definitions, $\div S_h\subset V_h$
and, using the commutativity, we have $\div S_h=V_h$, i.e.,
\raise.035in\hbox{$\begin{CD}S_h@>\div>> V_h \to 0\,\end{CD}$} is exact.  To complete this
to a discrete analogue of the elasticity complex, we define
$Y_h$ to be the inverse image of $S_h$ under $J$.  Then $Y_h$ is exactly
the space of $C^1$ piecewise quintic polynomials which are $C^2$ at the
vertices of the meshes.  This is in fact a well-known finite element
space, called the Hermite quintic or Argyris space, developed for solving
$4$th order partial differential equations (for which the inclusion in
$H^2(\Omega)$ and therefore $C^1$
continuity is required).  The shape functions are
$\mathbb P_5(T)$ and the $21$ degrees of freedom are the values of
the function and all its first and second partial derivatives at the
vertices and the integrals of the normal derivatives along edges.
We then have a \emph{discrete elasticity complex}
\begin{equation*}
\begin{CD}
\mathbb P_1\hookrightarrow\,@. Y_h @>J>> S_h @>\div>> V_h @.\to0,
\end{CD}
\end{equation*}
or, diagrammatically,
\begin{equation*}
\begin{CD}
\mathbb P_1\hookrightarrow\,@.
 \raise-.2in\hbox{\includegraphics[width=.65in]{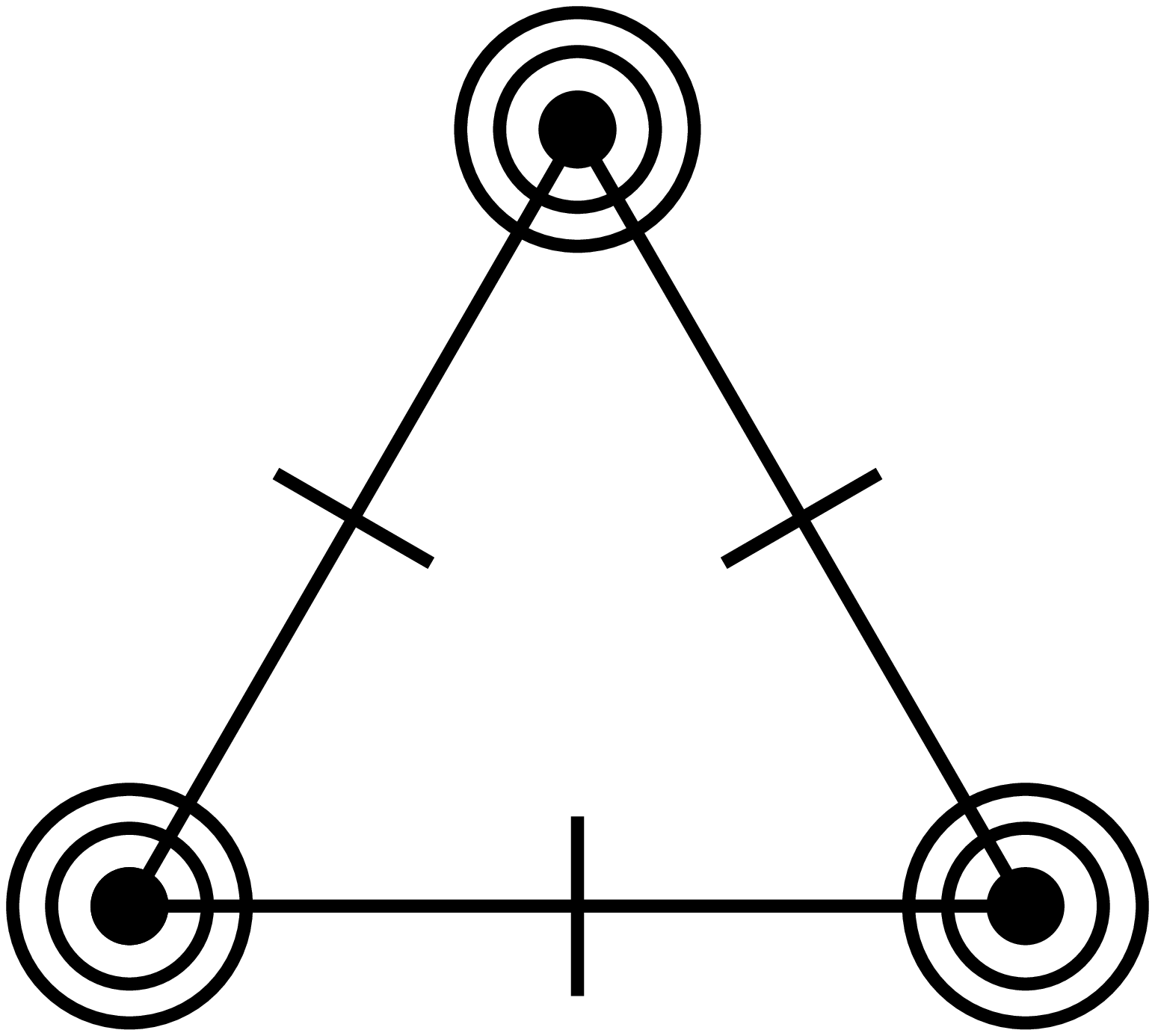}} @>J>>
 \raise-.19in\hbox{\includegraphics[width=.5in]{000233.eps}} @>\div>>
 \raise-.11in\hbox{\includegraphics[width=.5in]{000242.eps}} @.\to0.
\end{CD}
\end{equation*}
Moreover this sequence is exact and is coupled to the two-dimensional elasticity
sequence via a commuting diagram:
\begin{equation*}
\begin{CD}
\mathbb P_1\hookrightarrow\,@.
C^\infty(\Omega) @>J>> C^\infty(\Omega,\mathbb S) @>\div>>
C^\infty(\Omega,\R^3) @.\to0\\
 @.@VV\Pi^Y_hV @VV\Pi^S_hV @VV{\Pi^V_h}V\\
\mathbb P_1\hookrightarrow\,@. Y_h @>J>> S_h @>\div>> V_h @.\to0
\end{CD}
\end{equation*}
The right half of this diagram encodes the information necessary to
establish the stability of our mixed finite element method.

The Hermite quintic finite elements arose naturally from our mixed
finite elements to complete the commutative diagram.  Had they not
been long known, we could have used this procedure to devise
a finite element space contained in
$H^2(\Omega)$.  In fact, on close scrutiny we can see that any
stable mixed finite elements for elasticity with polynomial shape
functions will give rise to a finite element space with polynomial
shape functions contained in $H^2(\Omega)$.  However, it is known that
such spaces are difficult to construct and complicated.
In fact, it can be proved
that an $H^2$ finite element space must utilize shape functions of degree at least $5$
and the first and second partial derivatives at the
vertices must be among the degrees of freedom \cite{zenisek}.  This
helps explain why mixed finite elements for elasticity have proven
so hard to devise.  In particular, we can rigorously establish the
stress elements must involve polynomials of degree $3$, and that
vertex degrees of freedom are unavoidable.

In addition to the element just described, elements of all greater
orders are also introduced in \cite{mixedelas}.  The elements of next
higher order can be seen as the final two elements in this discrete
elasticity complex.
\begin{equation*}
\begin{CD}
\mathbb P_1\hookrightarrow\,@.
 \raise-.2in\hbox{\includegraphics[width=.65in]{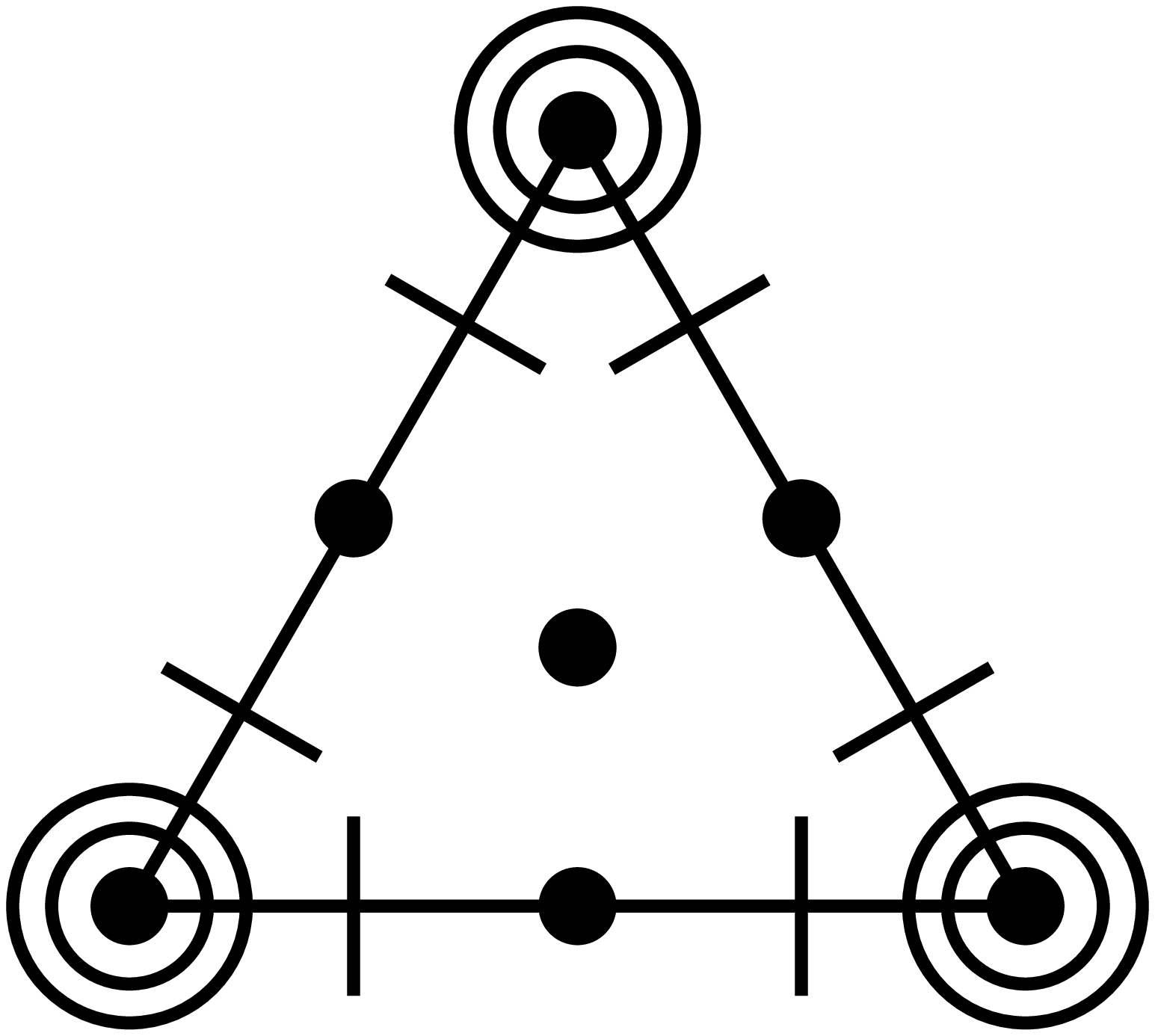}} @>J>>
 \raise-.19in\hbox{\includegraphics[width=.5in]{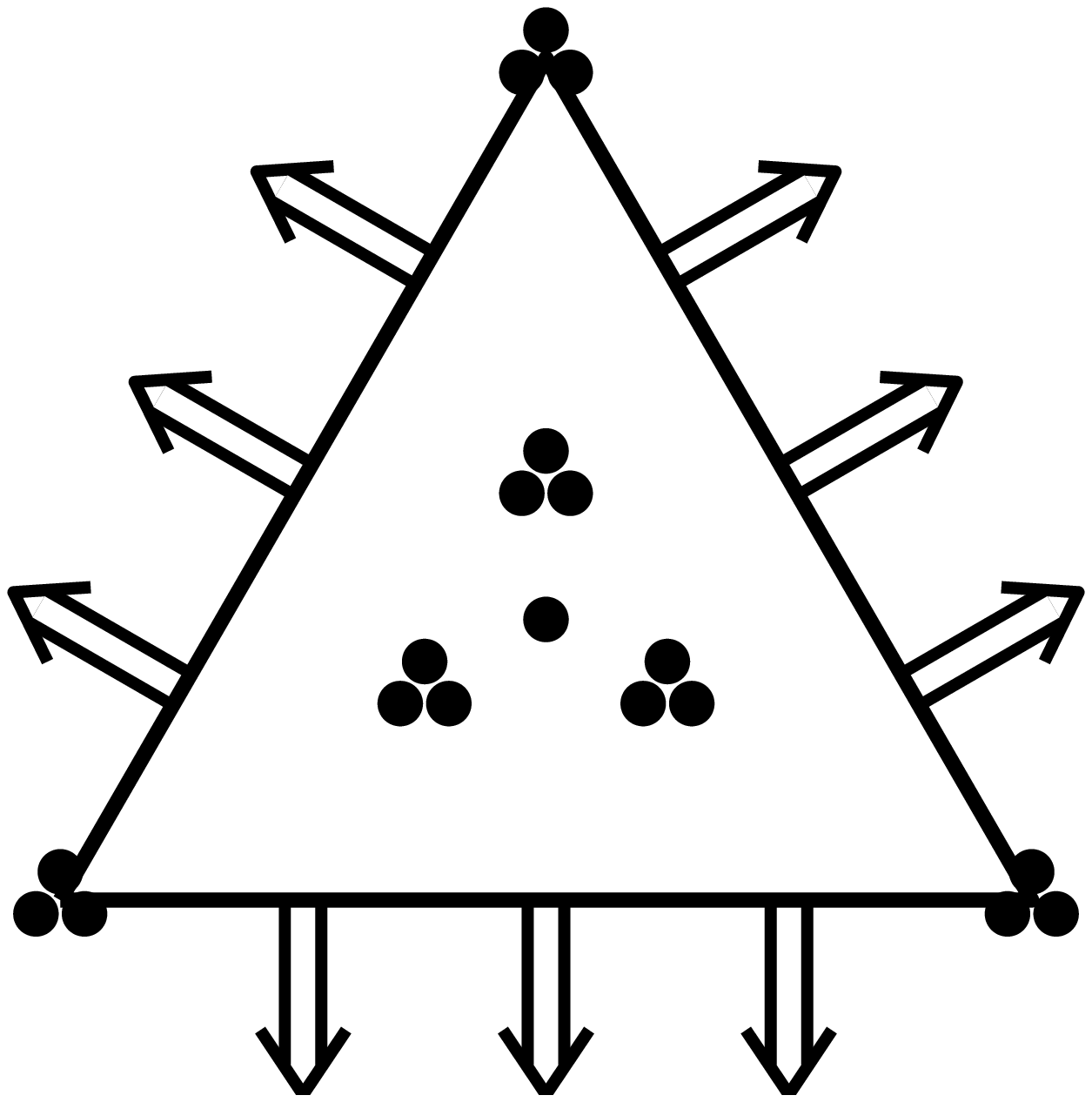}} @>\div>>
 \raise-.11in\hbox{\includegraphics[width=.5in]{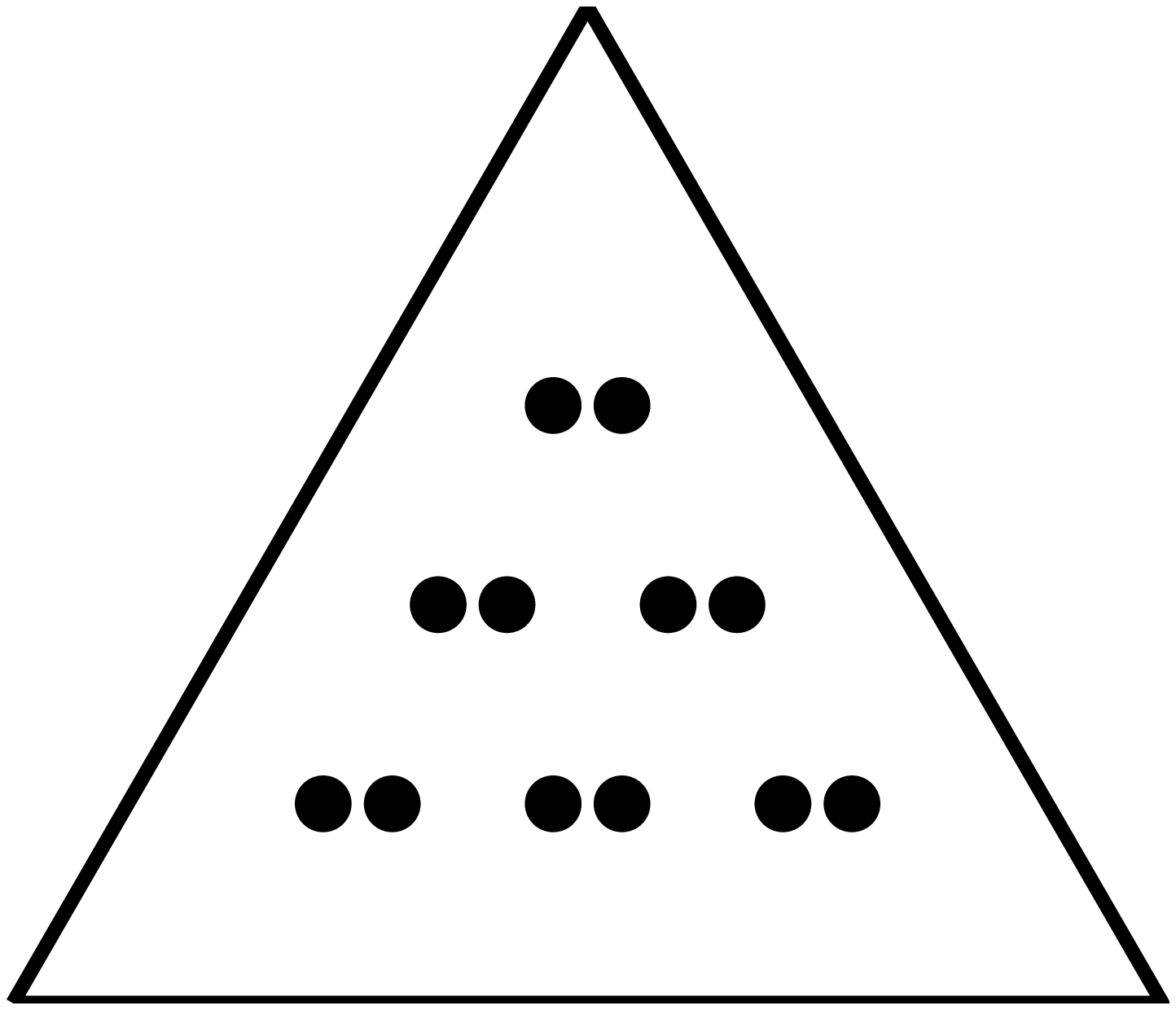}} @.\to0.
\end{CD}
\end{equation*}
It is also possible to simplify the lowest order element slightly.
To do this we reduce the displacement space from piecewise linear
vector fields to piecewise rigid motions, and we replace the
stress space with the inverse image under the divergence of the
reduced displacement space.  This leads to a stable
element shown in this exact sequence:
\begin{equation*}
\begin{CD}
\mathbb P_1\hookrightarrow\,@.
 \raise-.2in\hbox{\includegraphics[width=.65in]{000237.eps}} @>J>>
 \raise-.19in\hbox{\includegraphics[width=.5in]{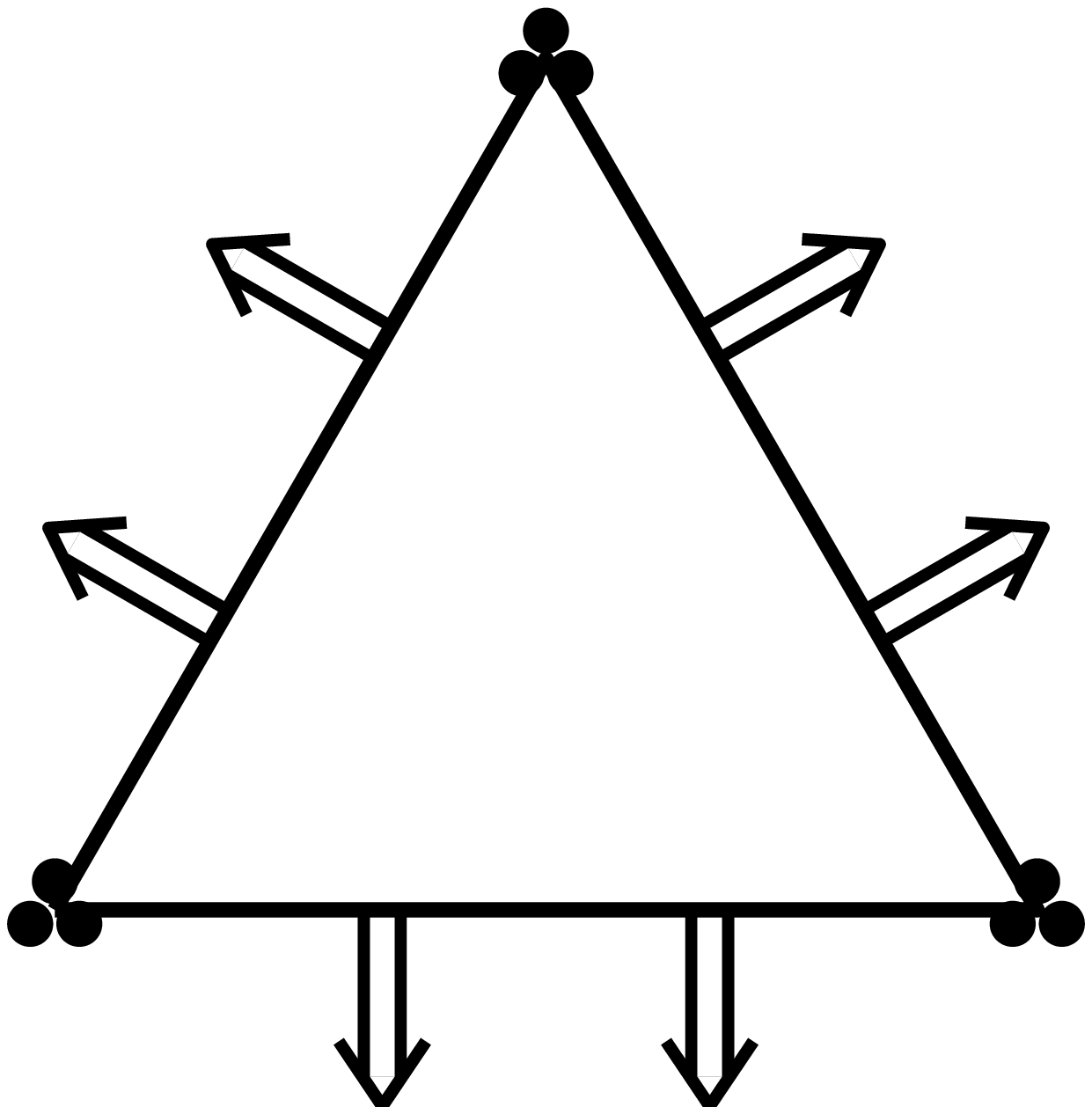}} @>\div>>
 \raise-.11in\hbox{\includegraphics[width=.5in]{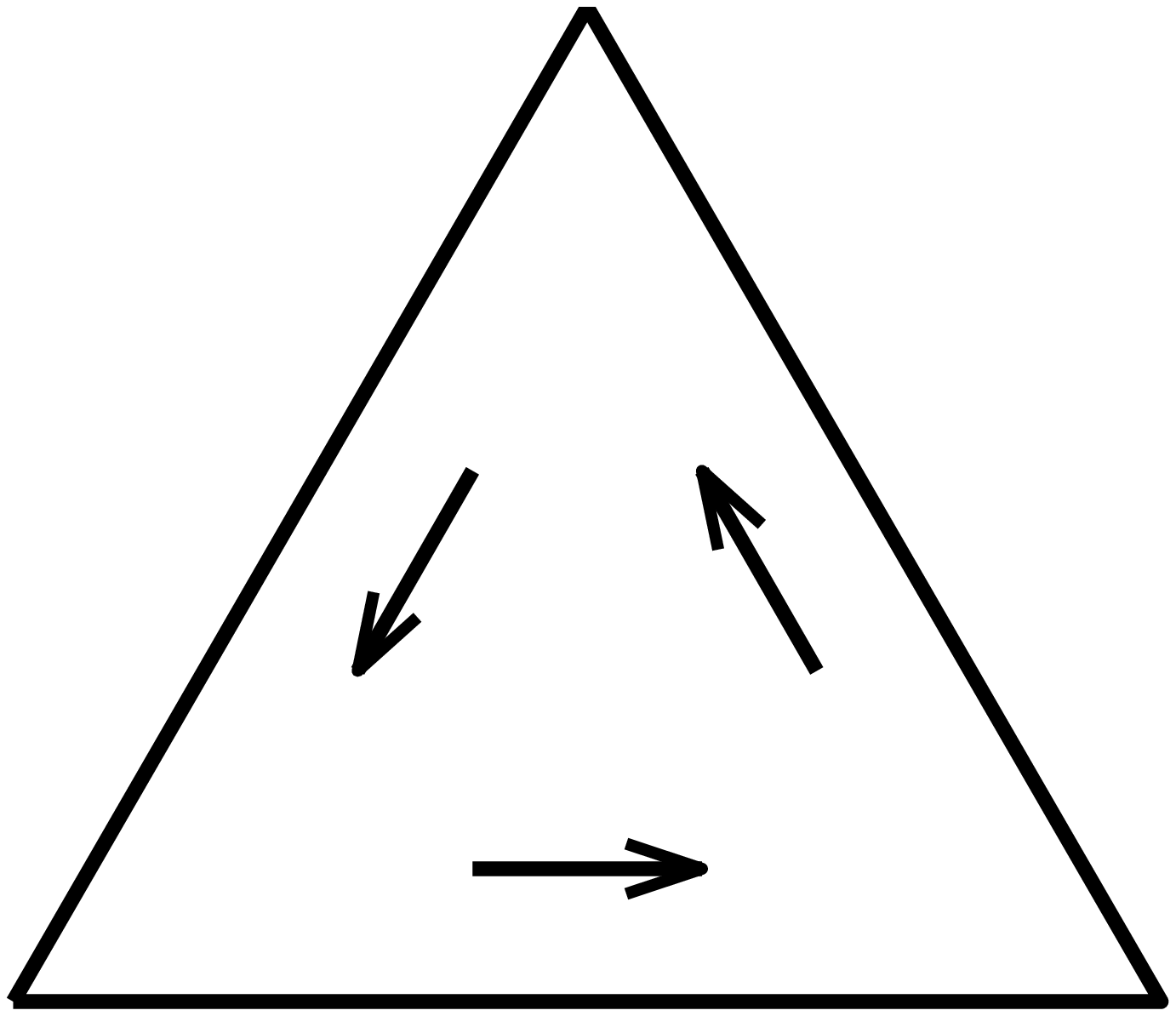}} @.\to0.
\end{CD}
\end{equation*}

Because of the unavoidable complexity of $H^2$ finite elements,
practitioners solving $4$th order equations often resort to
\emph{nonconforming} finite element approximations of $H^2$.  This
means that the finite element space does not belong to $H^2$ in that
the function or the normal derivative may jump across element
boundaries, but the spaces are designed so that jumps are small enough in some
sense (e.g., on average).  The error analysis is more complicated for
nonconforming elements, since in addition to stability and approximation
properties of the finite element space, one must analyze the
\emph{consistency error} arising from the jumps in the
finite elements.  In \cite{ncelas} Winther and the author investigated the
the possibility of nonconforming mixed finite elements for
elasticity, which, however are stable and convergent, and developed two
such elements.  These are
related to nonconforming $H^2$ elements via nonconforming discrete
elasticity complexes, two of which are pictured here:
\begin{equation*}
\begin{CD}
\mathbb P_1\hookrightarrow\,@.
 \raise-.2in\hbox{\includegraphics[width=.65in]{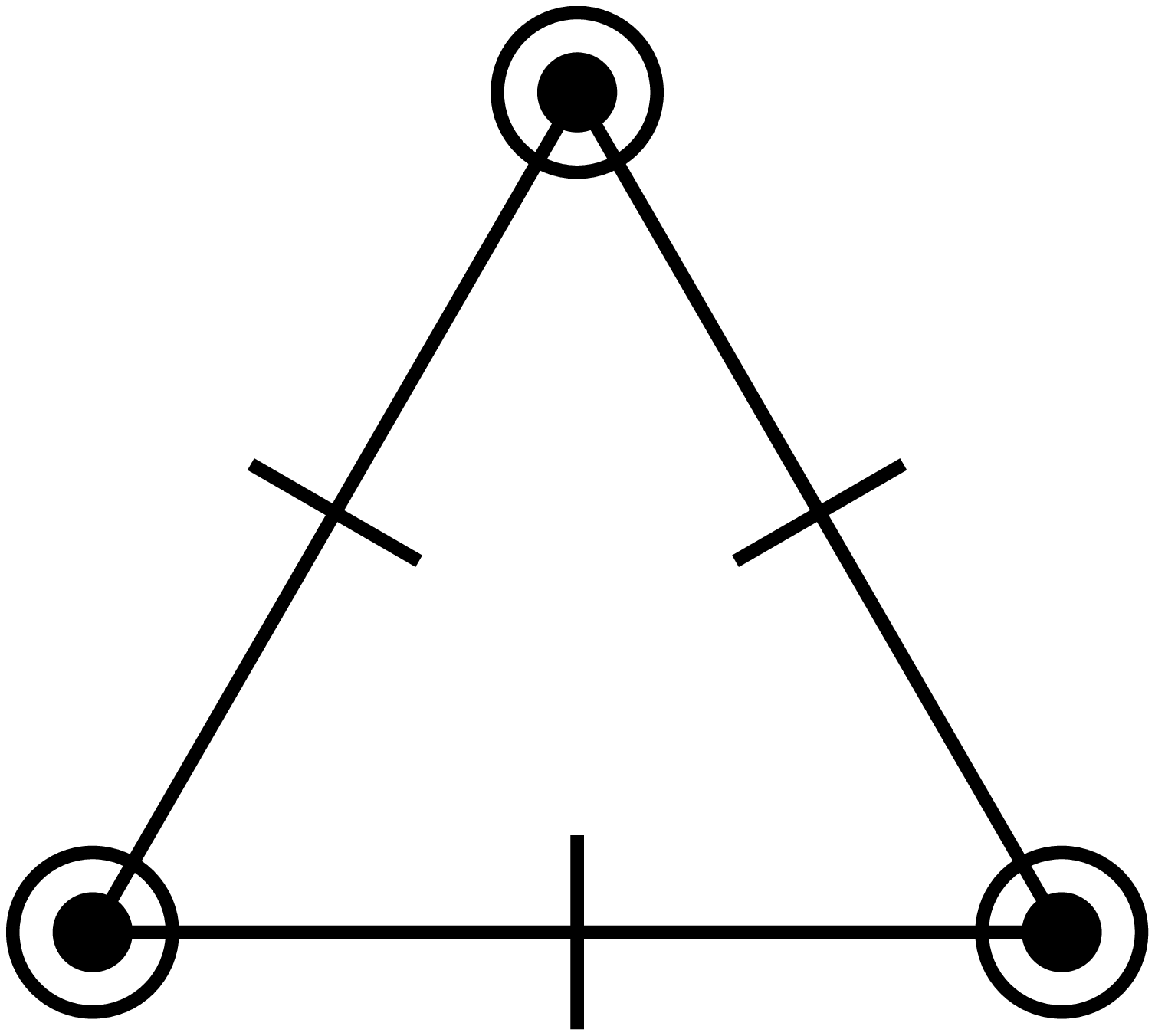}} @>J>>
 \raise-.22in\hbox{\includegraphics[width=.5in]{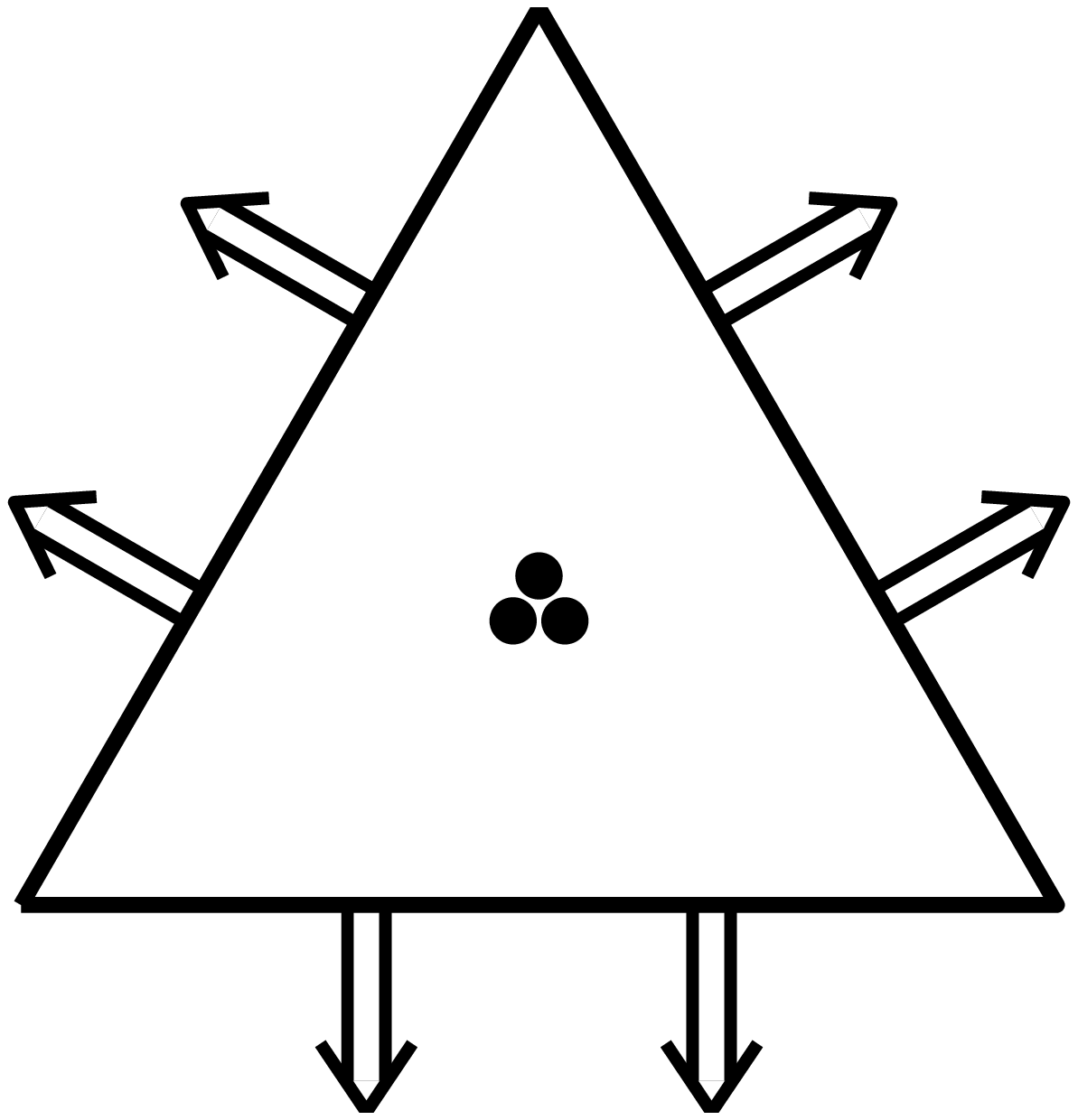}} @>\div>>
 \raise-.13in\hbox{\includegraphics[width=.5in]{000242.eps}} @.\to0.
\\
\\
\mathbb P_1\hookrightarrow\,@.
 \raise-.2in\hbox{\includegraphics[width=.65in]{000240.eps}} @>J>>
 \raise-.22in\hbox{\includegraphics[width=.5in]{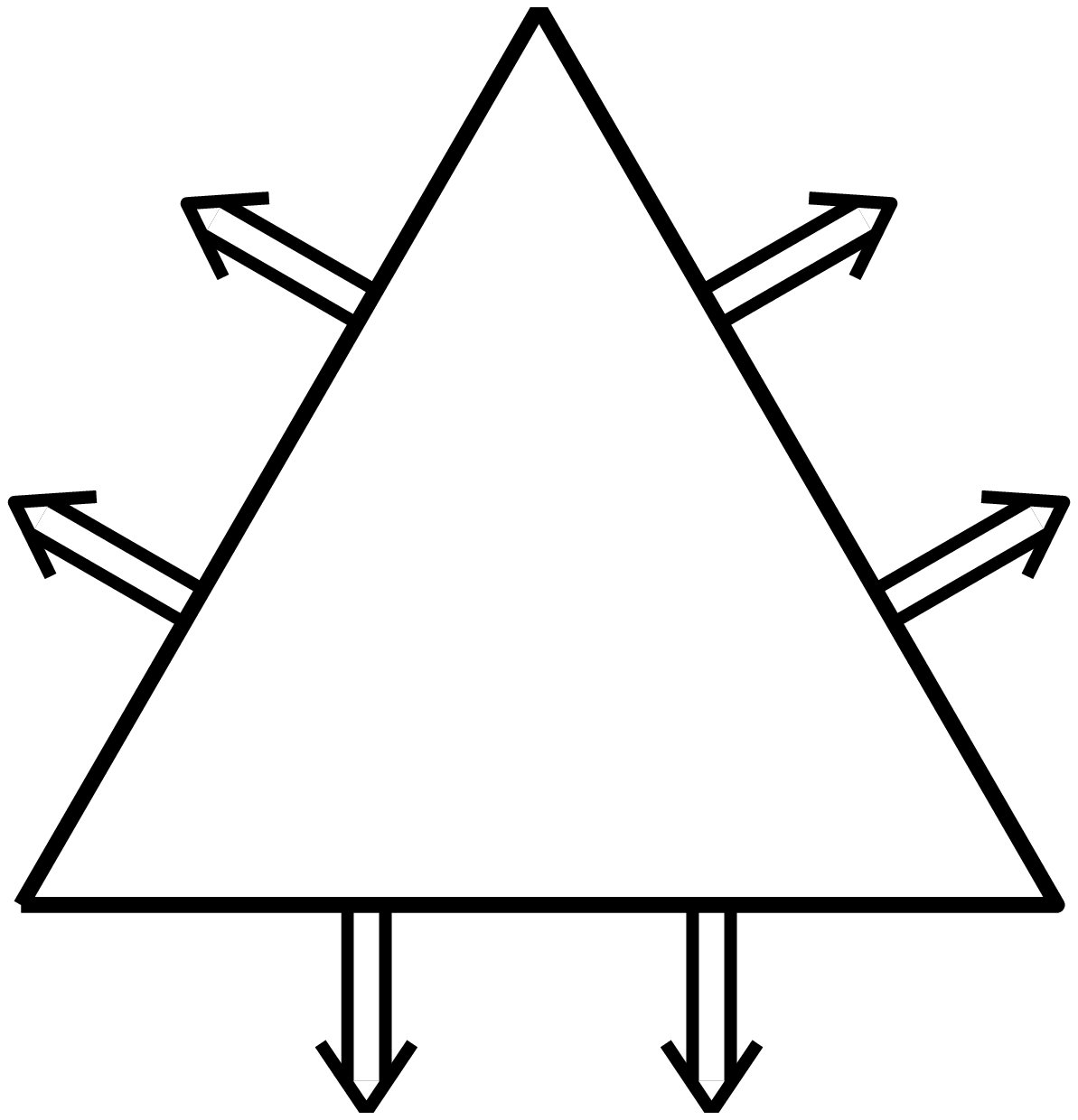}} @>\div>>
 \raise-.13in\hbox{\includegraphics[width=.5in]{000239.eps}} @.\to0.
\end{CD}
\end{equation*}
In both cases the shape function space for the stress is contained
between $\mathbb P_1(T,\mathbb S)$ and $\mathbb P_2(T,\mathbb S)$.
The nonconforming $H^2$ finite element depicted in these diagrams was
developed for certain $4$th order problems in \cite{ntw}.  Note
the nonconforming mixed elasticity elements are significantly
simpler than the conforming ones (and, in particular, don't require
vertex degrees of freedom).

\label{lastpage}

\end{document}